\documentclass[review]{elsarticle}

\usepackage[latin9]{inputenc}
\usepackage[english]{babel}
\usepackage{amstext,amsmath,amsfonts,amssymb,amsbsy,bm,amsthm}
\usepackage[usenames,dvipsnames]{color}
\usepackage{hyperref}
\usepackage{algorithmic}
\usepackage{ocr} % for the igatools' font
\usepackage{color}
\usepackage{mathtools}
\usepackage[lined,boxed,commentsnumbered]{algorithm2e}

\numberwithin{equation}{section}

\newtheorem{definition}{Definition}[section]
\newtheorem{remark}[definition]{Remark}

\newcommand{\R}{\mathbb{R}}

\newcommand{\ii}{\boldsymbol{i}}
\newcommand{\jj}{\boldsymbol{j}}
\newcommand{\kk}{\boldsymbol{k}}
\newcommand{\qq}{\boldsymbol{q}}
\newcommand{\xx}{\boldsymbol{x}}

\newcommand{\zzeta}{\boldsymbol{\zeta}}

\newcommand{\cchi}{\boldsymbol{\chi}}

\newcommand{\df}{\hat{\bm{D}} \bm{F}}

\def\B{\color{black}}

\usepackage{listings}

\lstnewenvironment{code}[1][]{%
   \lstset{%
      language=C++,
      basicstyle=\ttfamily\footnotesize,
      keywordstyle=\bfseries\color{blue},
      identifierstyle=\color{black},
      stringstyle=\color{red},
      commentstyle=\color{orange},
      morecomment=[l][\color{magenta}]{\#},
      captionpos=b,frame=single,
      numbers=left,                   
      escapeinside={\%*}{*)},
      #1,
   }
}{}

\begin{document}

\begin{frontmatter}

\title{Fast formation of isogeometric Galerkin matrices by weighted quadrature}

\author[cassino]{F. Calabr\`o} \ead{calabro@unicas.it}
\author[mat,imati]{G. Sangalli} \ead{giancarlo.sangalli@unipv.it}
\author[mat]{M. Tani} \ead{mattia.tani@unipv.it}

\address[cassino]{DIEI, Universit\`a degli Studi di Cassino e del Lazio Meridionale}

\address[mat]{Dipartimento di Matematica, Universit\`a degli Studi di
Pavia}

\address[imati]{Istituto di Matematica Applicata e Tecnologie
  Informatiche ``E. Magenes'' del CNR, Pavia}

\begin{abstract}
In this paper we propose an algorithm for the formation of matrices of
isogeometric Galerkin methods. The algorithm is based on three 
ideas. The first is that we perform the external loop over the
rows of the matrix. The second is that we calculate the row entries by 
weighted quadrature. The third is that we exploit the
(local) tensor product structure of the basis functions. While all
ingredients have a fundamental role for  computational efficiency, 
the major conceptual change of paradigm with respect to the standard
implementation is the idea of using weighted quadrature: the test function is incorporated in the integration
weight while 
the  trial function, the geometry parametrization and the PDEs
coefficients form the integrand function. This approach is very
effective in reducing the computational cost, while maintaining the
optimal order of approximation of the method. Analysis of the
cost is confirmed by numerical testing, where we show that, for $p$ large enough, the time
required by  the floating point operations is less than the time
spent in unavoidable memory operations (the sparse matrix allocation and memory
write). The proposed algorithm allows significant time saving when assembling isogeometric Galerkin matrices for all the 
degrees of the test spline space and  paves the way for a use of
high-degree $k$-refinement in isogeometric analysis. 
\end{abstract}


\end{frontmatter}

\section{Introduction}

 Isogeometric analysis has been introduced  by the seminal paper 
\cite{Hughes:2005} as    an extension of the classical
finite element method.  It is  based on the idea
of using the functions that are  adopted for the 
geometry parametrization in computer aided
design also    to represent  the numerical
solution of  the PDE of interest. These functions are   splines, Non-Uniform Rational B-Splines (NURBS) and
extensions. Many papers have demonstrated the advantage
of isogeometric methods  in various applications.  For the interested reader, we refer to  the book
\cite{Cottrell:2009}.

% A recent overview on the mathematical aspects of IGA is  \cite{acta},
% that covers the known mathematical theory of IGA but also contains
% an updated bibliography with references to the major contributions and applications of
% IGA in various engineering fields.

One interesting feature of isogeometric methods is the
possibility of using high-degree high-regularity splines as they  
deliver higher accuracy per degree-of-freedom in comparison to $C^0$
finite elements   \cite{Bazilevs:2006,BeiraoDaVeiga:2011,BeiraoDaVeiga:2014}. However, the computational cost per
degree-of-freedom is also higher for smooth splines, in currently
available  isogeometric codes.  In practice, quadratic or cubic
splines are  preferred as they  maximize 
computational efficiency. 

The computational cost of a  solver for a linear PDE  problem is the sum of
the cost of the formation of the system matrix and the cost of the
solution of the linear system. The former is dominant in standard 
isogeometric codes  already for  low degree  (see
e.g. \cite{Antolin:2015,sangalli2016isogeometric}). Recent papers in the literature 
have addressed this important issue (see
e.g. \cite{Mantzaflaris:2014,Rypl:2012}). 

In this paper we adopt the following definition of \emph{optimality}:
an algorithm for the formation of the matrix of a 
Galerkin method is \emph{optimal} if its computational cost is
of the order of  the number of non-zero entries of the matrix
to be calculated. Optimal algorithms  are known  in the case of $C^0$ finite elements
(see \cite{Ainsworth:2011,Ainsworth:2016}). However, this is
still an open problem for smooth splines.
 
We consider in this paper  a $d$-dimensional scalar Poisson model problem on a single-patch
domain, and an isogeometric tensor-product space of degree
$p$ and total dimension  $N_{\mathsf{DOF}}$, with $N_{\mathsf{DOF}}\gg
 p^d  $.  For the sake of simplicity, we focus on the case of
$C^{p-1}$ continuity,  i.e., the typical setting of the so-called $k$-method (see
e.g. \cite{Cottrell:2009}). The
resulting { stiffness matrix} has   $O( N_{\mathsf{DOF}}
(2p+1)^d) \approx   O( N_{\mathsf{DOF}} p^d)$ non-zero 
entries. Therefore, we assume $C N_{\mathsf{DOF}} p^d$ floating point\footnote{Throughout the
  paper, $C$ is a (reasonably small) constant that does not depend on  $N_{\mathsf{DOF}} $ and
$p$, and is in general  different at
each occurrence.} operations
(FLOPs) is the (quasi)-optimal 
computational cost for the formation of  the stiffness matrix.
The algorithms currently used in isogeometric codes are 
suboptimal with respect to the degree $p$, that is, their cost
grows with respect to  the degree $p$ faster than $p^d$.

The  majority of isogeometric codes inherit   a finite element
architecture, which adopt  an element-wise assembly loop with
element-wise standard  Gaussian quadrature (SGQ). 
Each local stiffness matrix has dimension $(p+1)^{2d}$
and each entry  is calculated by quadrature on $(p+1)^d$ Gauss points. The total cost is   $ O( N_{\mathsf{EL}} p^{3d})
\approx O( N_{\mathsf{DOF}} p^{3d})$  FLOPs, where $  N_{\mathsf{EL}} $  is the number of
elements and,  for the $k$-method, $  N_{\mathsf{EL}} \approx
N_{\mathsf{DOF}} $. 

The standard way  to reduce the cost  is to reduce the
number of quadrature points,  for example  by \emph{reduced Gaussian} 
\cite{Adam:2015b,Schillinger:2014}  (eventually corrected by
variationally consistent constraints \cite{hillman2015variationally})  or \emph{generalised Gaussian} quadrature (GGQ) \cite{Hughes:2010,bremer2010nonlinear,Cheng:1999,Ma:1996}. To clarify GGQ, consider  the  {mass matrix}
$\mathbb{M}= \{m_{\ii,\jj} \}  $ whose entries, calculated on the parametric domain ${\hat\Omega}= [0,1]^d$, have the form
\begin{equation} 
m_{\ii,\jj}=\int_{\hat\Omega} c(\zzeta)\, \hat{B}_{\ii} (\zzeta) \, \hat{B}_{\jj} (\zzeta) \,  d\zzeta   \ ,
\end{equation}
where $ \hat{B}_{\ii} $ and  $ \hat{B}_{\jj} $ are the tensor-product
B-spline,  and $c$  is a coefficient that incorporated
the determinant Jacobian of the geometry mapping and other possible
non-tensor product factors.  
The work \cite{Hughes:2010} has explored the 
possibility of constructing and using GGQ quadrature
 of the kind
\begin{equation} \label{eq:gG-quadrature}
\int_{\hat\Omega} c(\zzeta)\, \hat{B}_{\ii} (\zzeta) \hat{B}_{\jj}
(\zzeta) \,  d\zzeta   \approx  \mathbb{Q}^{\mathrm{GGQ}}(c(\cdot)\, \hat{B}_{\ii} (\cdot) \hat{B}_{\jj} (\cdot)),
\end{equation}
where the quadrature weights $w^{\mathrm{GGQ}}_{\qq}$  and  points $\xx^{\mathrm{GGQ}}_{\qq}$ of the
quadrature rule  $\mathbb{Q}^{\mathrm{GGQ}}( f (\cdot) ) = \sum_{\qq}
w^{\mathrm{GGQ}}_{\qq} f (\xx^{\mathrm{GGQ}}_{\qq})$ fulfil  the exactness conditions  
\begin{equation}\label{eq:gG-quadrature-exactness}
  \int_{\hat\Omega}  \hat{B}^{2p}_{\kk} (\zzeta) d\zzeta =
  \mathbb{Q}^{\mathrm{GGQ}}(\hat{B}^{2p}_{\kk} (\cdot) ), \quad \forall \kk.
\end{equation}
Here $\{ \hat{B}^{2p}_{\kk} (\cdot)  \}$ is the B-spline basis of
degree $2p$ and continuity $C^{p-1}$. Exact integration of product of
a pair of 
$p$ degree splines $\hat{B}_{\ii} (\cdot) \hat{B}_{\jj} (\cdot)$ is
then guaranteed by \eqref{eq:gG-quadrature-exactness}. Since the $w^{\mathrm{GGQ}}_{\qq}$
and  $\xx^{\mathrm{GGQ}}_{\qq}$ are not known analytically, they need to be
computed numerically as solution of the global non-linear problem
\eqref{eq:gG-quadrature-exactness}, see
\cite{AitHaddou:2015,Bartovn:2015b,bartovn2016optimal}, and the
recent paper \cite{Johannessen:2016} where the problem is 
effectively solved by a Newton method with continuation.
The paper \cite{Auricchio:2012} uses local exactness conditions
instead of \eqref{eq:gG-quadrature-exactness}.
 The number of  conditions in
 \eqref{eq:gG-quadrature-exactness} is  $\#\{ \hat{B}^{2p}_{\kk}
(\cdot)  \} \approx N_{\mathsf{DOF}} (p+1)^d \approx N_{\mathsf{EL}} (p+1)^d$, dropping the lower order
terms, and therefore  GGQ is expected to use  about $N_{\mathsf{EL}}
\left ( \frac{p+1}{2} \right ) ^d$ quadrature points, with a saving of
a factor $2^d$ with respect to SGQ.

The number of quadrature points is not the only issue to consider here,
and indeed  the element-wise assembling loop has  a relevant
role as well.   On one hand,  it  allows the reuse of finite
element available routines, which is a clear advantage as it greatly 
simplifies code development. On
the other hand, it is intrinsically  not optimal,  as each elemental
stiffness matrix  has size $(p+1)^d$ and
therefore the total computational cost is bounded from below by  $ C N_{\mathsf{EL}} p^{2d}
\approx C N_{\mathsf{DOF}} p^{2d}$.  This ideal threshold is
approached by \emph{sum-factorization}, that is, by  arranging
 the computations in a way that exploits the tensor-product
 structure of multivariate spline, with a cost of   $ O(
 N_{\mathsf{DOF}} p^{2d+1})$  FLOPS, see   \cite{Antolin:2015}.

Further  cost reduction is possible with a change of paradigm
from   element-wise  assembly. This has been explored in some recent
papers.  In \cite{Mantzaflaris:2015} the
integrand factor due to geometry and PDE coefficients is interpolated on the space of  splines shape functions on a uniform knot vector, the same space where the approximation is considered,  while the integrals arising are pre-computed in exact manner. The final cost of assembly in this case  is  $O( N_{\mathsf{DOF}}
p^{2d})$. Another  approach has been proposed in
\cite{mantzaflaris2014matrix} where the stiffness matrix is approximated
by a  low-rank sum of $R$ Kronecker matrices that can be efficiently formed
thanks to their structure. This is a  promising approach with
computational  cost of  $  O( N_{\mathsf{DOF}} R p^{d})$ FLOPs.

We propose in this paper a new  algorithm which does not use the
element-wise  assembling loop. Instead, we loop over the matrix rows
and we use a specifically designed    \emph{weighted  quadrature} (WQ) rule
for each row.
In particular, the quadrature rule for the  $\ii$-th row of
$\mathbb{M}$ is as follows:
\begin{equation} \label{eq:weighted-quadrature}
\int_{\hat\Omega} c(\zzeta)\, \hat{B}_{\jj} (\zzeta) \, (\hat{B}_{\ii}
(\zzeta)  d\zzeta )  \approx  \mathbb{Q}^{\mathrm{WQ}}_{\ii}(c(\cdot)\,
\hat{B}_{\jj} (\cdot)), \quad \forall  \jj.
\end{equation}
Unlike \eqref{eq:gG-quadrature}, in the right hand side of  \eqref{eq:weighted-quadrature}
the integrand function is $c(\cdot)\,
\hat{B}_{\jj} (\cdot)$ since the test function is incorporated into
the integral weight (measure) $(\hat{B}_{\ii} (\zzeta) d\zzeta) $.  The price to pay is that the
quadrature weights depend on $\ii$, while we select global quadrature points as
suitable interpolation points that do not depend  on $\ii$. Again, the
quadrature weights  are not known analytically and need to be
computed numerically as solution  of the
exactness conditions  
\begin{equation}\label{eq:weighted-quadrature-exactness}
  \int_{\hat\Omega} \, \, \hat{B}_{\jj} (\zzeta) \, (\hat{B}_{\ii} (\zzeta)
d\zzeta)  = \mathbb{Q}^{\mathrm{WQ}}_{\ii}(\hat{B}_{\jj}).
\end{equation}
However, the exactness conditions
\eqref{eq:weighted-quadrature-exactness} are linear with respect to the
weights. Furthermore, \eqref{eq:weighted-quadrature-exactness}  is a
local problem as the
weights outside $ \text{supp}( \hat{B}_{\ii} ) $ can be set to zero
a priori. The knot vectors do not need to be
uniform with this approach.

The number of exactness
conditions of  \eqref{eq:weighted-quadrature-exactness} is  $\#\{ \hat{B}_{\jj} (\cdot)  \} = N_{\mathsf{DOF}} 
$. This is lower than the number of conditions  of
\eqref{eq:gG-quadrature-exactness}, which is $\#\{ \hat{B}^{2p}_{\kk}
(\cdot)  \} \approx N_{\mathsf{DOF}} (p+1)^d$. Hence, the main advantage of
the WQ with respect to GGQ is that  that the former requires significantly fewer
quadrature  points. In the case of maximum regularity only $2$ points
are needed in each direction sufficiently far away from the boundary,
while $p+1$ points are taken on boundary knot-spans along directions
that end on the boundary. Adopting sum-factorization (see
\cite{Antolin:2015}),  the proposed  algorithm has a total
computational cost of  $  O( N_{\mathsf{DOF}} p^{d+1})$ FLOPs. 

In our numerical benchmarking, performed in MATLAB, we have compared the standard GeoPDEs
3.0 (see \cite{de2011geopdes,geopdesv3}) mass matrix formation, based
on  element-loop SGQ,  with our  row-loop WQ-based algorithm, showing the impressive advantage.
WQ speedup is more than a factor of four for quadratics and rapidly grows with
the degree $p$. For example, the mass matrix  on a $20^3$ elements grid is calculated in about 62 hours
vs 27 seconds for degree $p=10$. For high $p$,  the asymptotic growth of the
computational time  is slower than the estimated cost from FLOPs counting:  for all  degrees of
 practical interest  the  growth we have measured is at most  $C
 N_{\mathsf{DOF}} p^{d}$, that is, optimal. This is due to the fact that the memory operations
dominate: in particular we have verified that  the
matrix  formation time with our  implementation  is mainly used in allocation (MATLAB's \texttt{sparse} call) and
memory write at least for sufficiently high degree $p$.  
%There is little room for further improvements.

The  WQ we propose is designed in order to fit into the 
mathematical theory that guarantees optimal order of
convergence of the method. This theory  is based on the  Strang lemma
\cite{ciarlet2002finite,strang1973analysis}.   We do not enter into
this topic, which is technical, and postpone it to a further work. In
this paper we give numerical evidence of optimal convergence on a
simple 1D benchmark. 
% Strang lemma requires
% that quadrature error is of the same order of the approximation error
% when the solution  and the coefficient ($c(\cdot)$, in this case) are
% smooth enough. This condition is unsymmetric, as it is      condition \eqref{eq:weighted-quadrature-exatness} is the main ingredient in order to recover the optimal order
% of accuracy of the whole isogeometric method. 

We also do not investigate parallelisation  in this paper, however we think the
proposed algorithm is well suited for a parallel implementation
since each matrix row is calculated independently,  which should
alleviate the race condition of typical finite element  element-wise
assembly (see  \cite{Karatarakis:2014} for details).

The outline of the paper is as follows. 
In Section 2, we present the idea of the WQ rules for univariate B-splines.
In Section 3, we briefly discuss the use of isogeometric analysis on a model problem, and fix the notation for the following sections. 
In Section 4, we extend the construction of WQ rules to the multivariate case; a pseudo-code is presented in Section 5, where the computational cost is also discussed.
In Section 6, we give details on the application of the WQ rules for the formation of the mass matrix.
In Section 7 we test the procedure: a simple 1D
test is performed in order to confirm accuracy and tests are presented
in order to compare time needed for the formation of mass matrices in 3D.  
%In Section 2, with the aim of fixing notations, we provide a very simple introduction on the use of
%isogeometric analysis for a model problem. In Section 3 we discuss the
%construction of the WQ rules. In Section 4 we discuss their application for the formation of the mass matrix. In both these sections we provide pseudo-codes
%for the algorithms. In Section 5 we test the procedure: a simple 1D
%test is performed in order to confirm accuracy and tests are presented
%in order to compare time needed for the formation of mass matrices in  
%3D. 
A complete benchmarking of the proposed procedure is beyond
the scope of the present paper and will be the subject of a forthcoming
paper. Finally, in Section 8 we draw conclusions. 

\section{Weighted quadrature}\label{Sect:Weighted}
Assume we want to compute integrals of the kind:
\begin{eqnarray} \label{eq:int_1d_1}
\int_{0}^1 \hat{B}_{i} (\zeta) \, \hat{B}_{j} (\zeta) \,  d\zeta   ,
\end{eqnarray}
where $\{\hat{B}_{i}\}_{i=1,\dots , n_{\mathsf{DOF}} }$
are $p$-degree univariate B-spline basis
functions defined on the parametric patch $[0,1]$.
We denote by
$\cchi$ the knot vector of distinct knots that define the
univariate B-splines $ \hat  B_{i} (\zeta) $. Moreover we define
\emph{knot-spans}  as the intervals
$[\chi_{k},\chi_{k+1} ], k=1,\dots, n_{\mathsf{EL}}$,
where $n_{\mathsf{EL}}:=(\# \cchi)-1$. The knot vector 
\begin{equation}
\displaystyle \Xi := \left\{ \xi_1,\xi_2, \dots , \xi_{\scriptscriptstyle n_{\scriptscriptstyle \mathsf{KNT}}} \right\} \ .
  \label{Xi}  
  \end{equation}
that defines
the univariate B-splines contains knots with repetitions depending
on the  regularity: if a knot $\chi_{k}$ has
multiplicity  $p-r$ then the univariate spline is  $C^{r}$
continuous at $\chi_{k} $.  For simplicity, we consider $r=p-1$
throughout this paper. Though it is not difficult to consider
arbitrary $r$, the proposed strategy takes advantage of high regularity. In order to focus on the relevant properties, we restrict our attention in this section to the uniform knot-spans, i.e. $\chi_{k+1}-\chi_{k} =h$ $\forall k =1, \ldots,n_{\mathsf{EL}}-1$. Moreover, we do not consider boundary functions, so we assume that the knot vector is periodic. 
Being in the context of Galerkin method, $
\hat{B}_{i} (\zeta)$ is denoted as a  test function and $\hat{B}_{j}
(\zeta)$ as a trial function.

\begin{figure}
\begin{center}
 \includegraphics[width=0.49\textwidth] {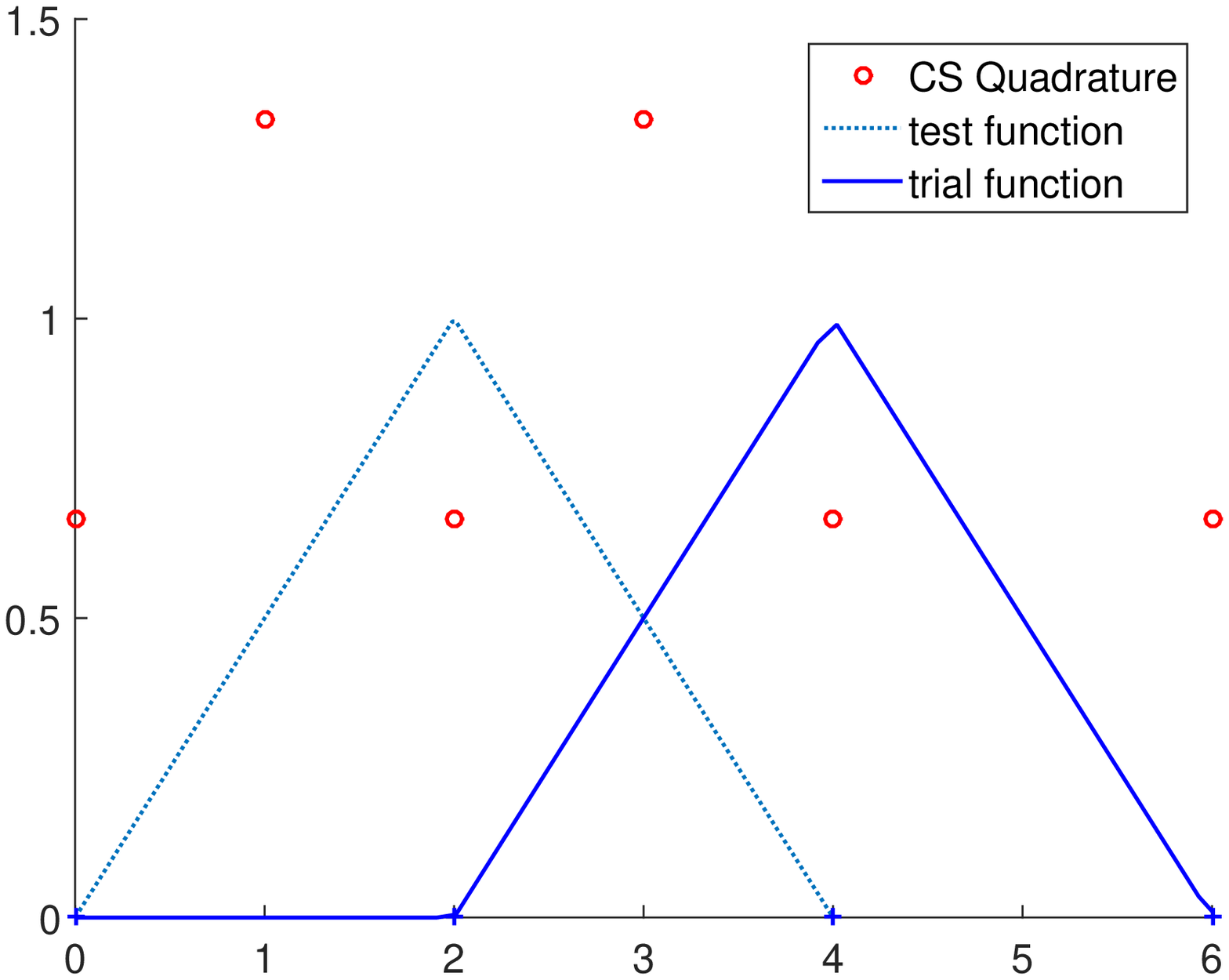} 
 \includegraphics[width=0.49\textwidth] {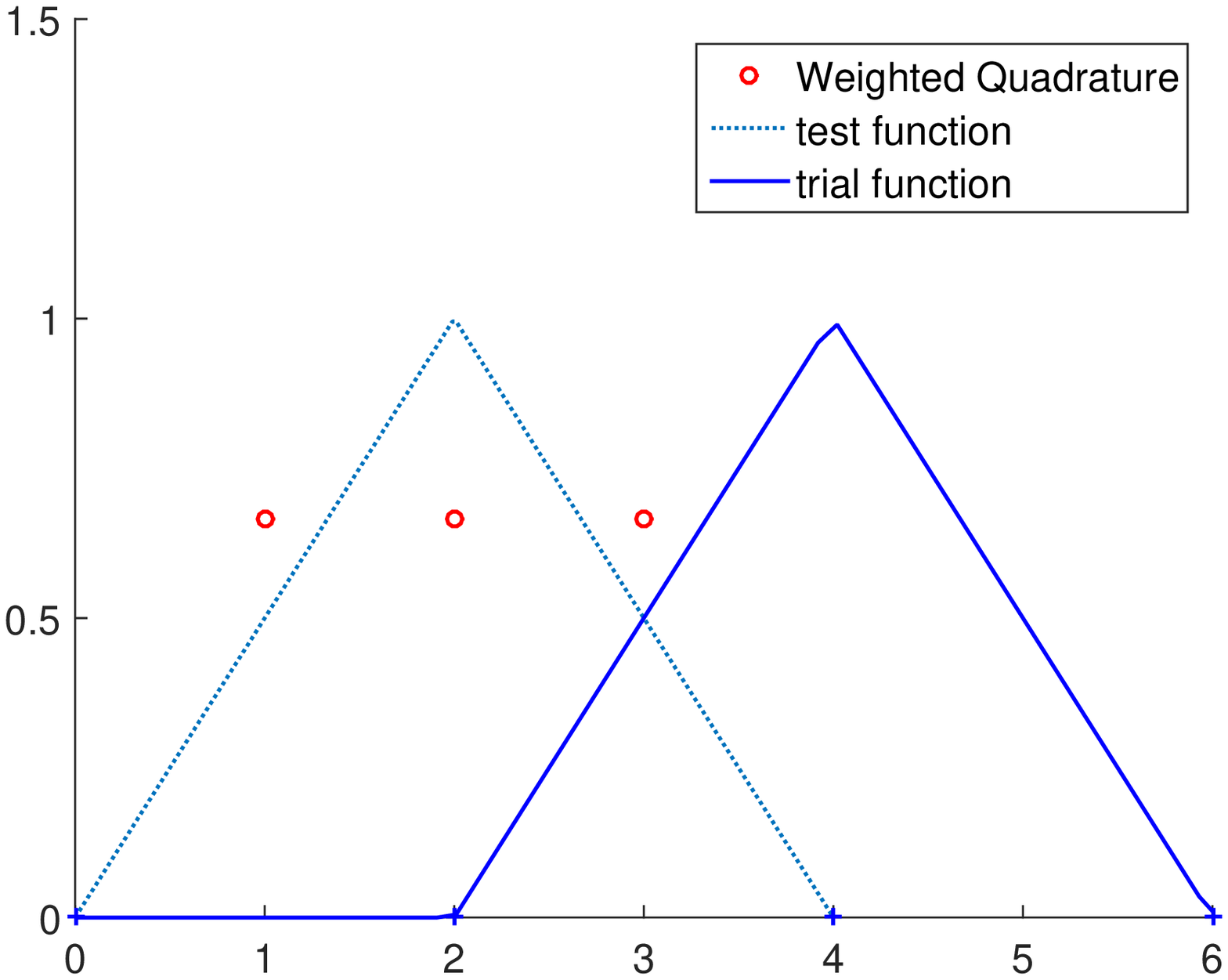} \\
\end{center}
\caption{Quadrature rule for B-spline with $p=1$. Cavalieri-Simpson quadrature rule (on the left) and it's interpretation as weighted quadrature (on the right). The active points and weights for $\mathbb{Q}_i^{WQ}$ are highlighted. In this and the next figure we set $h=2$ so that the quadrature points coincide with the integers, being the knots and the midpoints of the knot-spans.} \label{fig:WQ1d}
\end{figure}

We are interested in a fixed point quadrature rule. In the lowest
degree case, $p=1$, exact integration is performed by a composite
Cavalieri-Simpson rule (note that in this case this quadrature is also the Gauss-Lobatto $3$ points rule):
\begin{eqnarray} \label{eq:int_1d_2}
\int_{0}^1 \hat{B}_{i} (\zeta) \, \hat{B}_{j} (\zeta) \,  d\zeta  =  \mathbb{Q}^{CS} (\hat{B}_{i} \, \hat{B}_{j} ) = \sum_q w_q^{CS}\hat{B}_{i} (x_q^{CS}) \, \hat{B}_{j}(x_q^{CS}) ,
\end{eqnarray}
where $x_q^{CS}$ are the quadrature points and $w_q^{CS}$
the relative weights, see Figure \ref{fig:WQ1d}. In the above hypotheses the points $x_q^{CS}$ are the knots and the midpoints of the knot-spans and $w_q^{CS} = \frac{h}{3}$ on knots and $w_q^{CS} = \frac{2h}{3}$ on midpoints. 

Unbalancing the role of the test and the trial factors in \eqref{eq:int_1d_2}, we can see it as a weighted quadrature:
\begin{eqnarray} \label{eq:int_1d_3}
\int_{0}^1 \hat{B}_{i} (\zeta) \, \hat{B}_{j} (\zeta) \,  d\zeta  =  \mathbb{Q}_i^{WQ} (\hat{B}_{j} ) = \sum_q w_{q,i}^{WQ} \hat{B}_{j} (x_{q,i}^{WQ}) ,
\end{eqnarray}
where $x_q^{CS}=x_{q,i}^{WQ}$ and $w_{q,i}^{WQ}= \hat{B}_{i} (x_{q,i}^{WQ}) w_q^{CS}$. Because of the local support of the function $\hat{B}_{i}$ only in three points the quadrature $\mathbb{Q}_i^{WQ}$ is non-zero and the weights are equal to $\frac{h}{3}$ see Figure \ref{fig:WQ1d}.

If we go to higher degree, we need more quadrature points in
\eqref{eq:int_1d_2}. For $p$-degree splines the integrand $\hat{B}_{i}
\hat{B}_{j}$ is a piecewise polynomial of degree $2p$ and an
element-wise integration requires $2p+1$ equispaced points, or $p+1$
Gauss points, or about $p/2$ points with generalized Gaussian
integration (see
\cite{Hughes:2010,Auricchio:2012,bremer2010nonlinear,Calabro:2015}).
On the other hand, we can generalize \eqref{eq:int_1d_3} to higher
degree still using  as quadrature points only the knots and midpoints
of the knot spans.  
Indeed this choice ensures that, for each basis function
$\hat{B}_{i}$, $i=1,\ldots,n_{\mathsf{DOF}}$, there are $2p+1$
``active'' quadrature points where $\hat{B}_{i}$ is nonzero. Therefore
we can compute the  $2p+1$ quadrature weights by imposing
conditions for the  $2p+1$   B-splines $\hat{B}_{j}$ that need to be
exactly integrated.
Clearly, the advantage of the weighted quadrature approach is that its computational complexity, i.e., the total number of quadrature points, is independent of $p$.

%Once the quadrature points are selected, we can compute the weights from the exactness conditions. 
For the sake of clarity, we first consider the case $p=1$ in
detail. The exactness conditions  are:

%\begin{eqnarray}
\begin{flalign}
\label{eq:int_1d_4}
&\int_{0}^1 \hat{B}_{i} (\zeta) \, \hat{B}_{i-1} (\zeta) \,  d\zeta = \frac h 6 = \mathbb{Q}_i^{WQ} (\hat{B}_{i-1} ) = \frac 1 2 w_{1,i}^{WQ} , \\
&\int_{0}^1 \hat{B}_{i} ^2 (\zeta)\,  d\zeta = \frac{ 2 h}{ 3} = \mathbb{Q}_i^{WQ} (\hat{B}_{i} ) = \frac 1 2 w_{1,i}^{WQ} + w_{2,i}^{WQ} + \frac 1 2 w_{3,i}^{WQ} , \\
&\int_{0}^1 \hat{B}_{i} (\zeta) \, \hat{B}_{i+1} (\zeta) \,  d\zeta = \frac h 6 = \mathbb{Q}_i^{WQ} (\hat{B}_{i+1} ) = \frac 1 2 w_{3,i}^{WQ} \ .
\end{flalign}
%\end{eqnarray}
Then it is easy to compute $w_{q,i}^{WQ}= h/3 ,\, \forall q=1,2,3$. 

In the case $p=2$, five points are active and we have five exactness equations:
%\begin{eqnarray} 
\begin{flalign}
\label{eq:int_1d_5}
&\int_{0}^1 \hat{B}_{i} (\zeta) \, \hat{B}_{i-2} (\zeta) \,  d\zeta = \frac{ h}{120} = \mathbb{Q}_i^{WQ} (\hat{B}_{i-2} ) = \frac 1 8 w_{1,i}^{WQ} , \\
&\int_{0}^1 \hat{B}_{i} (\zeta) \, \hat{B}_{i-1} (\zeta) \,  d\zeta = \frac{26 h}{120} = \mathbb{Q}_i^{WQ} (\hat{B}_{i-1} ) = \frac 3 4 w_{1,i}^{WQ} + \frac 1 2 w_{2,i}^{WQ} + \frac 1 8 w_{3,i}^{WQ} , \\
&\int_{0}^1 \hat{B}_{i} ^2 (\zeta)\,  d\zeta = \frac{66 h}{120} = \mathbb{Q}_i^{WQ} (\hat{B}_{i} ) = \frac 1 8 w_{1,i}^{WQ} + \frac 1 2 w_{2,i}^{WQ} + \frac 3 4 w_{3,i}^{WQ} + \frac 1 2 w_{4,i}^{WQ} + \frac 1 8 w_{5,i}^{WQ} , \\
&\int_{0}^1 \hat{B}_{i} (\zeta) \, \hat{B}_{i+1} (\zeta) \,  d\zeta = \frac{26 h}{120} = \mathbb{Q}_i^{WQ} (\hat{B}_{i+1} ) = \frac 1 8 w_{3,i}^{WQ} + \frac 1 2 w_{4,i}^{WQ} + \frac 3 4 w_{5,i}^{WQ} , \\
&\int_{0}^1 \hat{B}_{i} (\zeta) \, \hat{B}_{i+2} (\zeta) \,  d\zeta = \frac{h}{120} = \mathbb{Q}_i^{WQ} (\hat{B}_{i+2} ) = \frac 1 8 w_{5,i}^{WQ} \ .
%\end{eqnarray}
\end{flalign}
In the previous calculation we have used the usual properties of B-splines that can be found, e.g., in \cite[Section 4.4]{Schumaker:2007}.
The system can be solved and leads to the following solution
$w_{q,i}^{WQ}= \frac{h}{30} \left[ 2, 7, 12, 7, 2 \right] $. 

\begin{figure}
\begin{center}
 \includegraphics[width=0.49\textwidth] {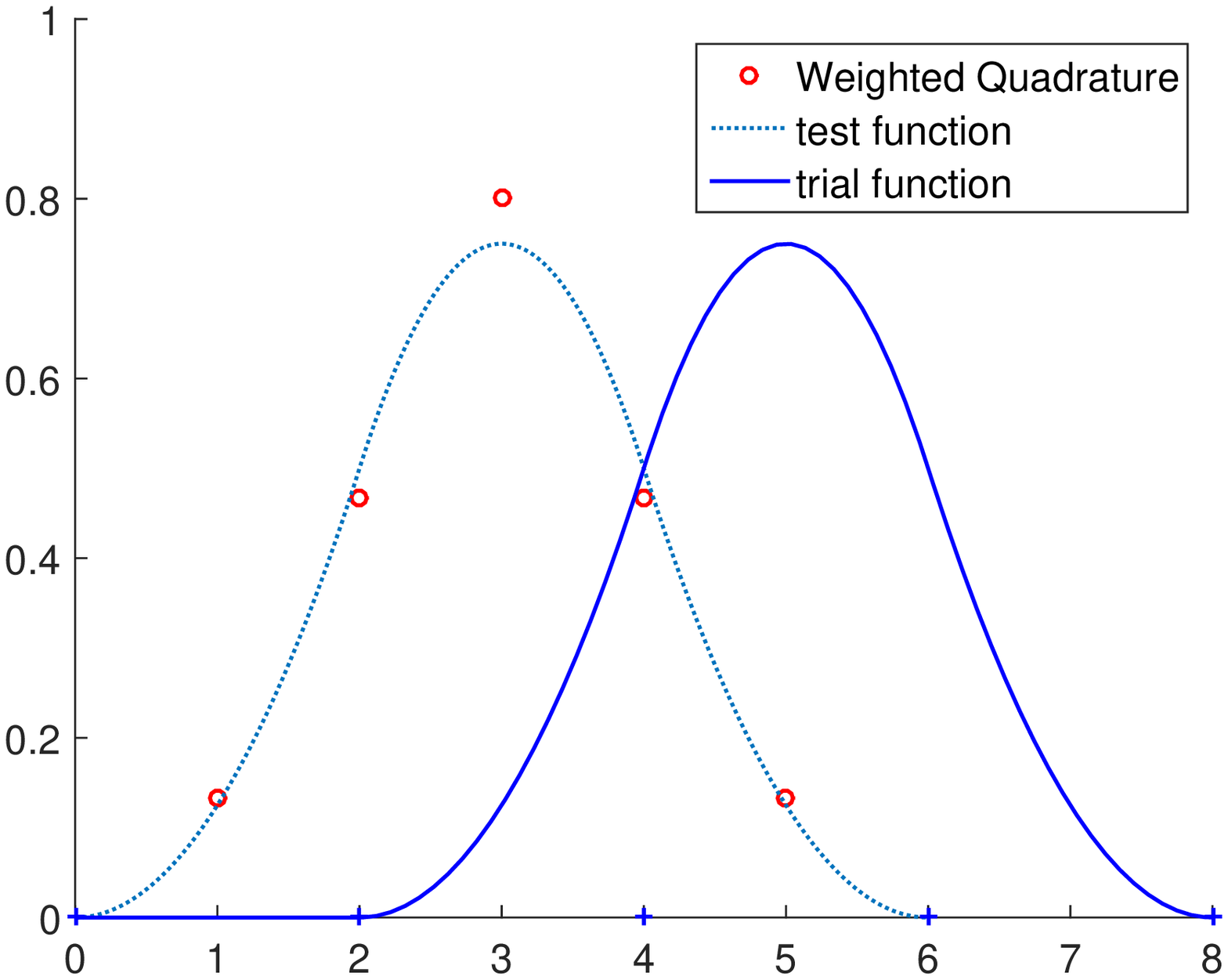} 
 \includegraphics[width=0.49\textwidth] {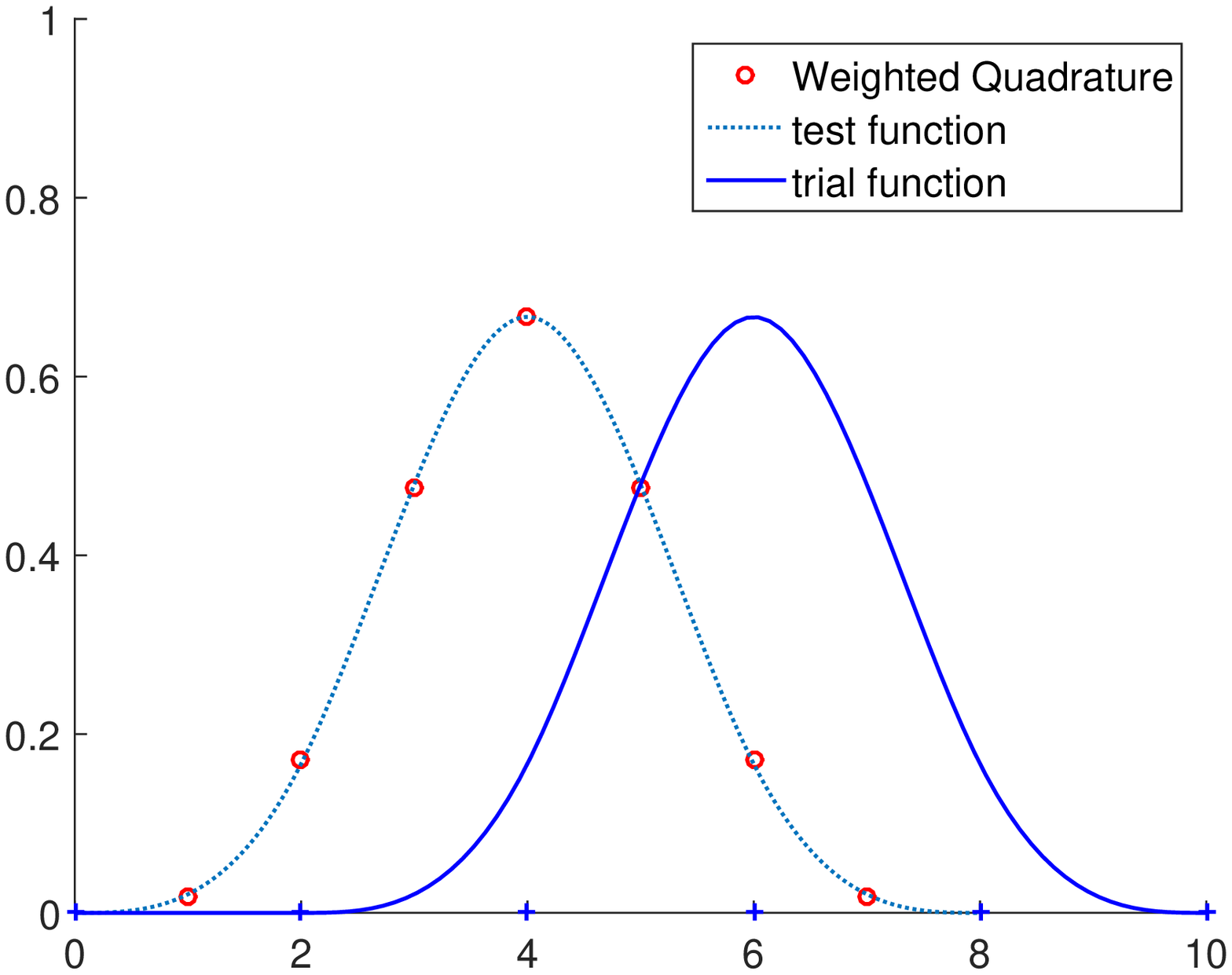} \\
 \includegraphics[width=0.49\textwidth] {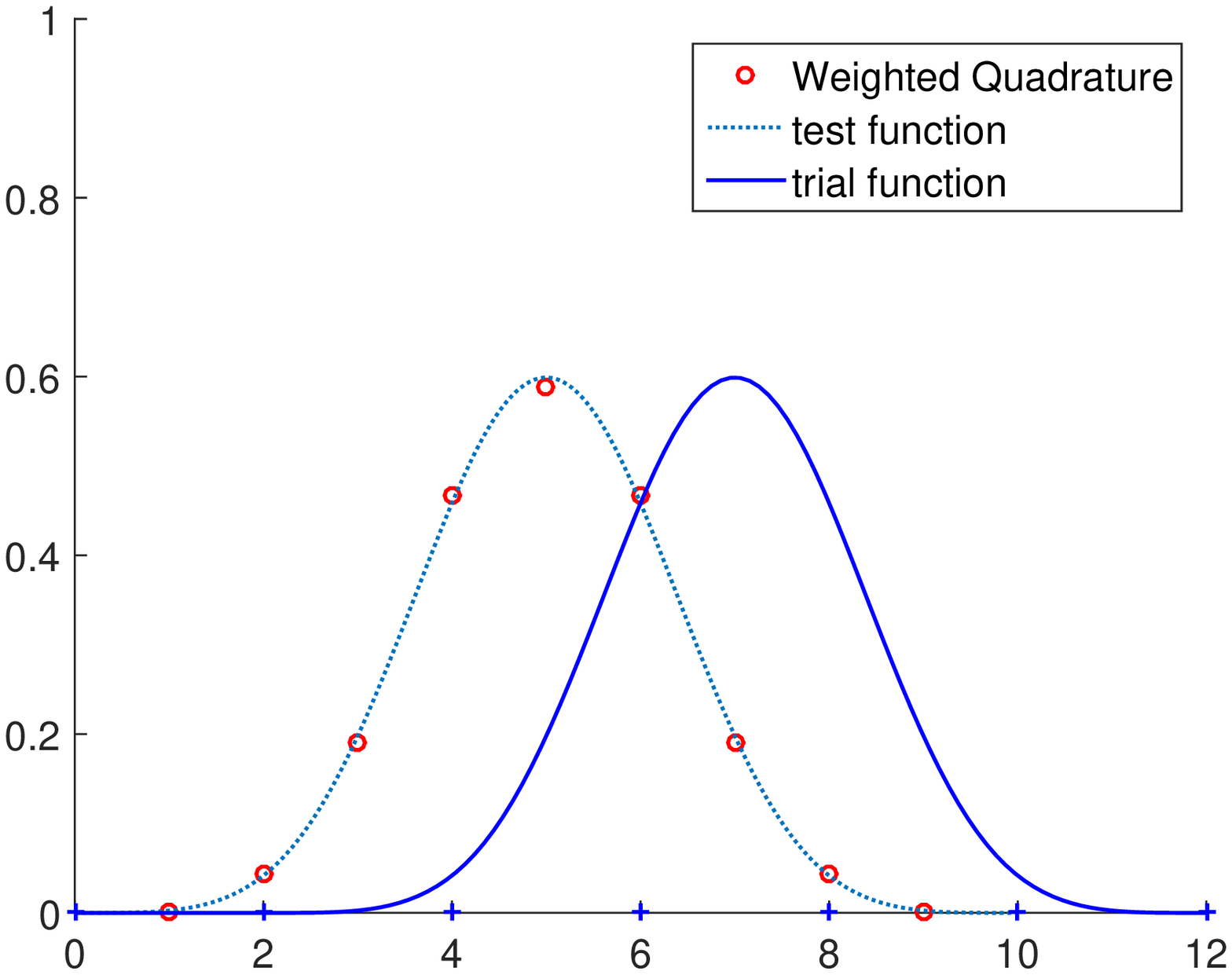} 
 \includegraphics[width=0.49\textwidth] {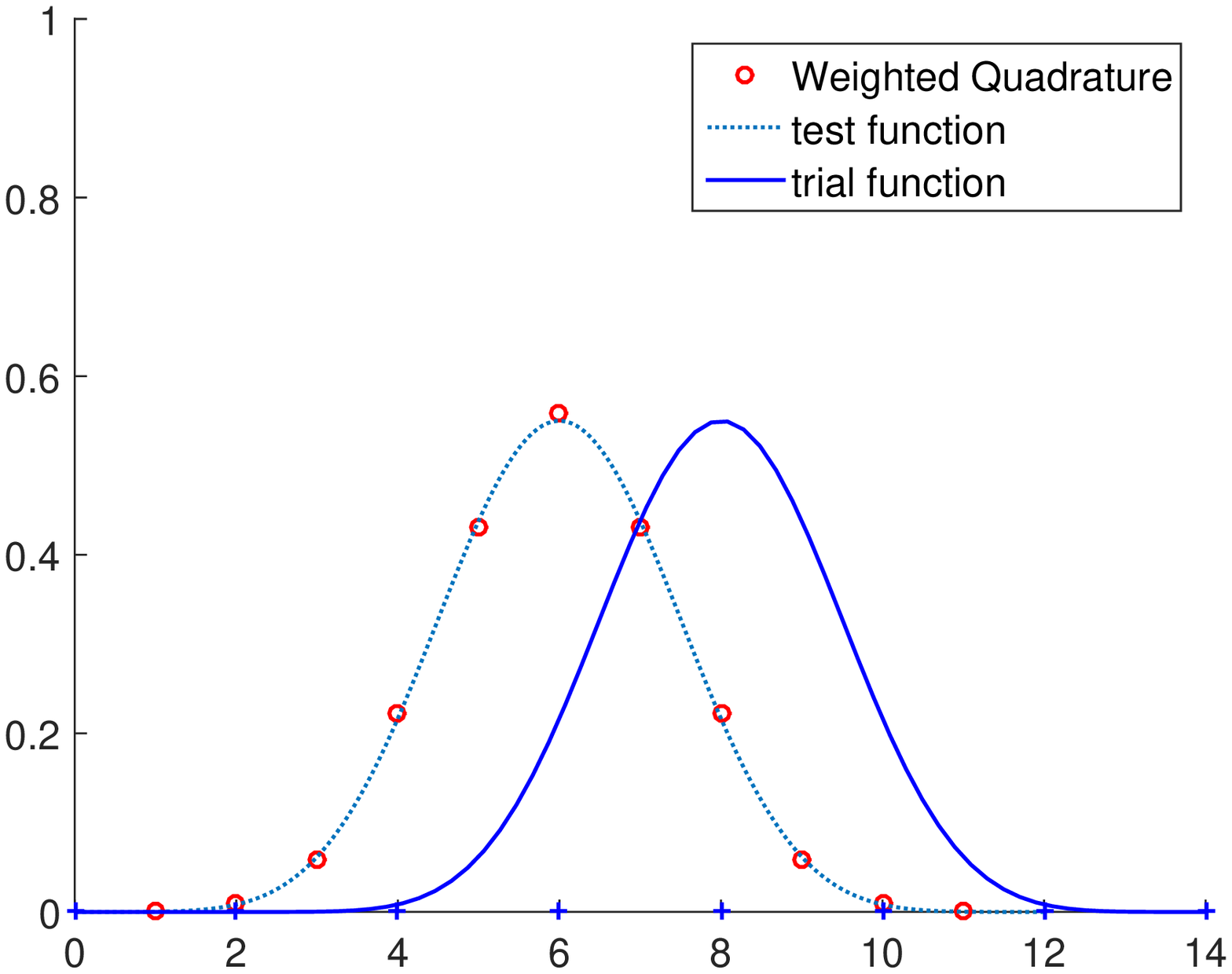} \\
\end{center}
\caption{Weighted quadrature rule for B-spline with various degrees, namely $p=2$ (upper left), $p=3$ (upper right), $p=4$ (lower left) and $p=5$ (lower right). 
%The active points and weights for $\mathbb{Q}_i^{WQ}$ are hightlighted. On the top line, on the left the case $p=2$, on the right the case $p=3$, on the bottom line on the left the case $p=4$, on the right the case $p=5$. 
Interestingly, we see that the weights displace around the values of the basis function (up to the scale factor $h/2$, which in this case is simply 1). In particular, this gives numerical evidence of the positivity of the weights, which in turn implies the stability of the rules.
%The displacement of weights around the values of the basis function becomes evident and gives the empirical measure of stability of the rules, thou not a formal proof. 
} \label{fig:WQ1d_2}
\end{figure}

When $p=3$, the same approach  gives $w_{q,i}^{WQ}= {h} \left[
  \frac{1}{105} , \frac{3}{35}, \frac{5}{21},
  \frac{1}{3},\frac{5}{21}, \frac{3}{35}, \frac{1}{105} \right]$.
These computations are reported in Figure \ref{fig:WQ1d_2}. 

In general case (arbitrary degree and non-uniform  spacing, boundary functions, lower regularity ...) the rule can be computed numerically as solution of a linear system, see Section \ref{Sect:Pseudocodes}.

%\begin{figure}
%\begin{center}
% \includegraphics[width=0.45\textwidth] {WQp4.eps} 
% \includegraphics[width=0.45\textwidth] {WQp5.eps} \\
%\end{center}
%\caption{Weighted quadrature rule for B-spline with $p=4,5$. The active points and weights for $\mathbb{Q}_i^{WQ}$ are hightlighted: on the left the case $p=4$, on the right the case $p=5$.} \label{fig:WQ1d_2b}
%\end{figure}

Given a weighted quadrature rule of the kind above, we are then interested in using it for the approximate calculation of integrals as:

\begin{eqnarray} \label{eq:int_1d_b1}
\int_{0}^1 c (\zeta) \hat{B}_{i} (\zeta) \, \hat{B}_{j} (\zeta) \,  d\zeta   \approx \mathbb{Q}_i^{WQ} \left(  c(\cdot) \hat{B}_{j}(\cdot) \right)= \sum_q w_{q,i}^{WQ} c(x_{q,i}^{WQ}) \hat{B}_{j} (x_{q,i}^{WQ}) \ .
\end{eqnarray}
For a non-constant function $c(\cdot)$, \eqref{eq:int_1d_b1} is in general just an approximation.
In particular, the symmetry of the integral is not preserved, that is $\mathbb{Q}_i^{WQ} \left(  c(\cdot) \hat{B}_{j}(\cdot) \right)$ is different from $\mathbb{Q}_j^{WQ} \left(  c(\cdot) \hat{B}_{i}(\cdot) \right)$. Consider, for example, the case $p=3$ derived above and apply the weighted quadrature rules to the linear function $c (\zeta)= \zeta$ in the case $j=i+1$. For simplicity we take $h=2$ so that the quadrature points are $x_{q,i}^{WQ}=[1:7]$. Then:
\begin{eqnarray*} \label{eq:int_1d_b2}
\mathbb{Q}_i^{WQ} \left(  c(\cdot) \hat{B}_{j}(\cdot) \right) = %\frac{2}{105} 1 0 + \frac{6}{35} 2 0 + 
\frac{10}{21} 3 \frac{1}{48} +  \frac{2}{3} 4 \frac{1}{6} + \frac{10}{21} 5 \frac{23}{48} + \frac{6}{35} 6 \frac{2}{3} + \frac{2}{105} 7 \frac{23}{48} \approx  2.3647, \\
\mathbb{Q}_j^{WQ} \left(  c(\cdot) \hat{B}_{i}(\cdot) \right) = \frac{2}{105} 3 \frac{23}{48} + \frac{6}{35} 4 \frac{2}{3} + \frac{10}{21} 5 \frac{23}{48} +  \frac{2}{3} 6 \frac{1}{6} + \frac{10}{21} 7  \frac{1}{48} %+ \frac{6}{35} 8 0 + \frac{2}{105} 9 0 
\approx  2.3615.
\end{eqnarray*}

A detailed mathematical analysis of the quadrature error 
%and symmetry error 
of weighted quadrature is of key interest, especially in the context of isogeometric Galerkin methods. This is however beyond the scope of this paper and for its importance deserves future work.

\section{Integral arising in isogeometric Galerkin methods}

We consider the model reaction-diffusion problem
\begin{equation}\label{poisson}
\left\{
\begin{aligned}
-\nabla^2 u + u = f&\quad \text{on} \qquad \Omega , \\
u = 0 & \quad \text{on} \qquad \partial \Omega ,
\end{aligned}
\right.
\end{equation}
Its Galerkin approximation on a discrete space $V$
requires the formation of the stiffness matrix $S$ and mass matrix $M$ whose entries are 
\begin{equation}
  \label{eq:stiffness-3d}
 s_{\ii,\jj} =  \int_{\Omega} \nabla R_{\boldsymbol{i}}(\boldsymbol{x})
  \nabla R_{\boldsymbol{j}}(\boldsymbol{x}) d\boldsymbol{x},
\end{equation}

\begin{equation}
  \label{eq:mass-3d}
  m_{\ii,\jj} = \int_{\Omega}  R_{\boldsymbol{i}}(\boldsymbol{x})
   R_{\boldsymbol{j}}(\boldsymbol{x}) d\boldsymbol{x},
\end{equation}

$R_{\boldsymbol{i}}$ and $R_{\boldsymbol{j}}$ being two basis functions in $V$. The dimension of the space $V$ is $N_{\mathsf{DOF}}:= \#V$.
Non-constant coefficients could be included as well.

In the isogeometric framework, $\Omega$ is given by a spline or NURBS
parametrization. Our notation follows  \cite[Section
4]{BeiraoDaVeiga:2014} and \cite{Antolin:2015}. For the sake of simplicity, we assume $\Omega$ is
given by a $d$-dimensional single patch spline representation, then it
is of the form:
\begin{displaymath}
  \Omega = \boldsymbol{F} (\hat\Omega), \text{ with }
  \boldsymbol{F}(\boldsymbol{\zeta}) = \sum_{\boldsymbol{i}}
  \boldsymbol{C}_{\boldsymbol{i}} \hat  B_{\boldsymbol{i}} (\boldsymbol{\zeta}),
\end{displaymath}
where $  \boldsymbol{C}_{\boldsymbol{i}} $ are the control points and
$  \hat B_{\boldsymbol{i}}$

%Copiato sopra, modifica per il multi-d are $p$-degree multivariate B-spline basis
% functions defined on the parametric patch $\hat\Omega = [0,1]^d$.

 We denote by
$\cchi_{l}$ the knot vector of distinct knots that define the
univariate B-splines $ \hat  B_{i_l} (\zeta_l) $
along the $l$-th direction. For each direction we have
\emph{knot-spans}  as the intervals
$[\chi_{l,k},\chi_{l,k+1} ], k=1,\dots, n_{\mathsf{EL},l}$,
where $n_{\mathsf{EL},l}:=(\# \cchi_{l})-1$. By cartesian
product,  they form a mesh of
$N_{\mathsf{EL}}= \prod_{l=1}^d n_{\mathsf{EL},l}$ \emph{elements} on
$\hat\Omega$. 
%We also define for each direction the knot vector that defines the univariate B-splines,
%The knot vector 
The knot vector, with possibly repeated knots, that defines the univariate B-spline space in the $l-$th direction is denoted as
\begin{equation}
\displaystyle \Xi_{l} := \left\{ \xi_1,\xi_2, \dots , \xi_{\scriptscriptstyle n_{\scriptscriptstyle \mathsf{KNT},l}} \right\} \ .
  \label{Xi_md},  
  \end{equation}
% that defines the univariate B-splines 
% contains knots with repetitions depending on the  regularity: if a knot $\chi_{k}$ has
% multiplicity  $p-r$ then the univariate splines is  $C^{r}$
% continuous at $\chi_{k} $.  For simplicity, we consider $r=p-1$
% throughout this paper. Though it is not difficult to consider
% arbitrary $r$, the proposed strategy takes advantage of high regularity. \B
As in Section \ref{Sect:Weighted}, we restrict to the case of maximum
regularity and allow repeated knots only at the endpoints of the open
knot vector.

The number of knots in each knot vector $n_{\mathsf{KNT},l}:=\# \Xi_{l}$
  is related to the number of degrees of freedom by
  $n_{\mathsf{DOF},l} +(p+1)= n_{\mathsf{KNT},l} $.  No assumptions are made on the length of the elements.   In our
  FLOPs counts, we always assume $p \ll n_{\mathsf{DOF},l}$, and then
  $n_{\mathsf{DOF},l}\approx n_{\mathsf{KNT},l} \approx
  n_{\mathsf{EL},l}$.

The multivariate B-splines are tensor-product of univariate B-splines:
\begin{equation}
  \label{eq:multivariate-B-splines}
  \hat B_{\ii}(\boldsymbol \zeta) = \hat B_{i_1}(\zeta_1) \ldots \hat B_{i_d}(\zeta_d).
\end{equation}

Above,  $\ii=(i_1,\ldots, i_d)$ is a multi-index that, with abuse of
notation, is occasionally as a scalar index, as in $
\ii=1,\ldots,N_{\mathsf{DOF}} $, with the relation $\ii \equiv
1+ \sum_{l = 1}^{d} (n_{\mathsf{DOF},1}\ldots
n_{\mathsf{DOF},l-1})(i_{l}-1)$.  We have  $N_{\mathsf{DOF}}=\prod_{l=1}^d n_{\mathsf{DOF},l}$.

Based on the
isogeometric/isoparametric paradigm, the basis functions $R_{\boldsymbol{i}}$ used in \eqref{eq:stiffness-3d}-\eqref{eq:mass-3d} 
are defined as $R_{\boldsymbol{i}} = \hat B_{\boldsymbol{i}}\circ
\boldsymbol{F} ^{-1}$; integrals are computed by change
of variable.  Summarizing, we are interested in the computation of \eqref{eq:mass-3d} after change of variable,   $\mathbb{M}= \{m_{\ii,\jj} \}      \in \mathbb{R}^{N_{\mathsf{DOF}} \times N_{\mathsf{DOF}} }$ where: 
\begin{eqnarray*}
 \ m_{\ii,\jj}=\int_{\hat\Omega} \hat{B}_{\ii} \, \hat{B}_{\jj}\, \text{det} \df \,  d\zzeta  \ .
\end{eqnarray*}
For notational convenience we write:  
\begin{eqnarray} \label{eq:mass2}
m_{\ii,\jj}=\int_{\hat\Omega} \hat{B}_{\ii} (\zzeta) \, \hat{B}_{\jj} (\zzeta) \, c(\zzeta)\,  d\zzeta   \ .
\end{eqnarray}
In more general cases, the factor  $c$ incorporates the  coefficient of the
equation and, for NURBS functions, the  polynomial denominator.
Similarly for the stiffness matrix $\mathbb{S}= \{s_{\ii,\jj} \}\in \mathbb{R}^{N_{\mathsf{DOF}} \times N_{\mathsf{DOF}} }$ we have:
\begin{eqnarray*}
s_{\ii,\jj} & = & \int_{\hat\Omega} \left( \df^{-T}
  \hat{\nabla} \hat{B}_{\ii}\right)^T  \left( \df^{-T}  \hat{\nabla} \hat{B}_{\jj} \right)  \text{det} \df \,  d\zzeta \\
& = & \int_{\hat\Omega} \hat{\nabla} \hat{B}_{\ii} ^T 
\left( \bigl[ \df^{-1} \df^{-T} \bigr] \text{det} \df \, \right)  \hat{\nabla} \hat{B}_{\jj} \,  d\zzeta \ 
\end{eqnarray*}
which we  write in compact form:
\begin{eqnarray} \label{eq:stiff}
  \begin{aligned}
    s_{\ii,\jj}&= \sum_{l,m=1}^d \int_{\hat\Omega} \left(\hat{\nabla}
      \hat{B}_{\ii} (\zzeta) \right)_l   c_{l,m}(\zzeta) \left(
      \hat{\nabla} \hat{B}_{\jj} (\zzeta) \right)_m  \,  d\zzeta .
  \end{aligned}
\end{eqnarray}
Here we have denoted by  $ \bigl\{ c_{l,m} (\zzeta)
\bigr\}_{ l,m=1,\dots ,d} $ the following matrix: 
\begin{eqnarray}\label{eq:coef_stiff}
c_{l,m} (\zzeta) = \bigl\{ \bigl[ \df^{-1}(\zzeta) \df^{-T}(\zzeta) \bigr] \,  \text{det} \df  (\zzeta) \bigr\}_{ l,m} .
\end{eqnarray}

%and expanding the matrix-vector products, as
%\begin{eqnarray} \label{eq:stiff}
%  \begin{aligned}
%    s_{\ii,\jj}&= \sum_{l,m=1}^d \int_{\hat\Omega} \left(\hat{\nabla}
%      \hat{B}_{\ii} (\zzeta) \right)_l   c_{l,m}(\zzeta) \left(
%      \hat{\nabla} \hat{B}_{\jj} (\zzeta) \right)_m  \,  d\zzeta \\
%&=\int_{\hat\Omega}  \hat{B}^\prime_{i_1} (\zeta_1) \hat{B}_{i_2}
%(\zeta_2)\cdots \hat{B}_{i_d} (\zeta_d)  c_{1,1}(\zzeta)
%\hat{B}^\prime_{j_1} (\zeta_1) \hat{B}_{j_2} (\zeta_2)\cdots
%\hat{B}_{j_d} (\zeta_d)  \,  d\zzeta   \\  
%&\quad + \int_{\hat\Omega}  \hat{B}_{i_1} (\zeta_1) \hat{B}^\prime_{i_2}
%(\zeta_2)\cdots \hat{B}_{i_d} (\zeta_d)  c_{2,1}(\zzeta)
%\hat{B}^\prime_{j_1} (\zeta_1) \hat{B}_{j_2} (\zeta_2)\cdots
%\hat{B}_{j_d} (\zeta_d)  \,  d\zzeta    \\ 
% &\quad + \ldots \\
% &\quad + \int_{\hat\Omega}  \hat{B}_{i_1} (\zeta_1) \hat{B}_{i_2}
%(\zeta_2)\cdots \hat{B}^\prime_{i_d} (\zeta_d)  c_{d,1}(\zzeta)
%\hat{B}^\prime_{j_1} (\zeta_1) \hat{B}_{j_2} (\zeta_2)\cdots
%\hat{B}_{j_d} (\zeta_d)  \,  d\zzeta  \\
% &\quad + \ldots  + \ldots \\
% &\quad + \int_{\hat\Omega}  \hat{B}_{i_1} (\zeta_1) \hat{B}_{i_2}
%(\zeta_2)\cdots \hat{B}^\prime_{i_d} (\zeta_d)  c_{d,d}(\zzeta)
%\hat{B}_{j_1} (\zeta_1) \hat{B}_{j_2} (\zeta_2)\cdots
%\hat{B}^\prime_{j_d} (\zeta_d)  \,  d\zzeta.
%  \end{aligned}
%\end{eqnarray}

The number of non-zero elements $N_{\mathsf{NZ}}$ of  $\mathbb{M}$ and
$ \mathbb{S}$ (the same for simplicity) depends on
the polynomial degree $p$ and the required regularity $r$.  
We introduce the following sets, where the support is
considered an open set:
\begin{gather}
\mathcal{K}_{l,i_l} =\left\{k\in \{1,\dots,n_{ \mathsf{EL},l}\} \,
  s.t.\,  ]\chi_{k-1}, \chi_{k} [  \subset \mbox{supp}\left(\hat{B}_{i_l} \right) \right\}\,,\\   
\mathcal{I}_{l,i_l} =\left\{j_l \in \{1,\dots,n_{ \mathsf{DOF},l}\}
\,   s.t.\,  \hat{B}_{i_l}\cdot \hat{B}_{j_l} \neq 0  \right\} \ ;
\end{gather}
and the related multi-indexes as:
\begin{gather}
\mathcal{K}_{\ii} = \prod_{l=1}^{d} \mathcal{K}_{l,i_l}\, \qquad \mathcal{I}_{\ii} =  \prod_{l=1}^{d} \mathcal{I}_{l,i_l} \ .
\end{gather}
We have  $ \# \mathcal{I}_{l,i} \le (2p+1) $ and  $ N_{\mathsf{NZ}}
= O( N_{\mathsf{DOF}}\, p^d) $. In particular, with maximal regularity in the case $d=1$ one has $N_{\mathsf{NZ}} = (2p+1) N_{\mathsf{DOF}}-p(p+1)$.
% When counting operations, we assume that the sets $\mathcal{K}_{\ii} ,\mathcal{I}_{\ii} $ are  provided by the meshing procedure, although our experience is that the cost for the construction of such connectivity information are negligible with respect to the other operations.

\section{Computation of the WQ rules} \label{sec:weighquad}

Consider the calculation of the mass matrix \eqref{eq:mass-3d}. The first step is to write the integral in a nested way, as done in \cite{Antolin:2015}:

\begin{multline}
m_{\ii,\jj}=\int_{{\hat\Omega}} \hat{B}_{\ii} (\zzeta) \hat{B}_{\jj} (\zzeta) c(\zzeta) \, d\zzeta = \\ \notag
\int_{0}^{1} \hat{B}_{i_1} (\zeta_1) \hat{B}_{j_1} (\zeta_1) 
\left[ \int_{0}^{1} \hat{B}_{i_2} (\zeta_2) \hat{B}_{j_2} (\zeta_2) \cdots  
\left[ \int_{0}^{1} \hat{B}_{i_d} (\zeta_d) \hat{B}_{j_d} (\zeta_d) c(\zzeta) \, d \zeta_d 
\right] \cdots  d \zeta_2\right] d \zeta_1  \notag
\end{multline}

Our idea is to isolate the \textit{test function} $\hat{B}_{i_l}$
univariate factors in each univariate integral  and to consider
it as a weight for the construction of the weighted  quadrature (WQ)
rule. This leads to a quadrature rule for each $i_l$  that is:
\begin{equation}\label{eq:Mass_approx}
  \begin{aligned}
    m_{\ii, \jj}\approx \tilde m_{\ii, \jj} & =
    \mathbb{Q}^{WQ}_{\ii}\left(    \hat{B}_{\jj }(\zeta)
      c(\zzeta)\right)    = \mathbb{Q}_{\ii}
    \left(    \hat{B}_{\jj }(\zeta) c(\zzeta)\right) \\
& = \mathbb{Q}_{i_1} \left( \hat{B}_{j_1}(\zeta_{1}) \mathbb{Q}_{i_2} \left( \cdots \mathbb{Q}_{i_d} \left(  \hat{B}_{j_d}(\zeta_{d}) c(\zzeta)  \right)   \right) \right) \  .
  \end{aligned}
\end{equation}
 Notice
that  we drop from now on the  label WQ used in the introduction in
order to simplify notation. The
key ingredients for the construction of the quadrature rules that
preserve the optimal approximation properties are the exactness
requirements. Roughly speaking,  exactness means that in
~\eqref{eq:Mass_approx} we have $m_{\ii, \jj}= \tilde m_{\ii, \jj}  $
whenever $c$ is a constant coefficient. When the stiffness term is
considered, also terms with derivatives have to be considered.

 We introduce the notation: 
\begin{equation} \label{eq:integrals}
\displaystyle  \begin{array}{l}
\displaystyle \mathbb{I}^{(0,0)}_{l,i_l,j_l} := \int_0^1 \hat{B}_{i_l}(\zeta_l)\hat{B}_{j_l}(\zeta_l)\, d\zeta_l  \\
\displaystyle \mathbb{I}^{(1,0)}_{l,i_l,j_l} := \int_0^1 \hat{B}^{\prime}_{i_l}(\zeta_l)\hat{B}_{j_l}(\zeta_l)\, d\zeta_l  \\
\displaystyle \mathbb{I}^{(0,1)}_{l,i_l,j_l} := \int_0^1 \hat{B}_{i_l}(\zeta_l)\hat{B}^{\prime}_{j_l}(\zeta_l)\, d\zeta_l  \\
\displaystyle \mathbb{I}^{(1,1)}_{l,i_l,j_l} := \int_0^1 \hat{B}^{\prime}_{i_l}(\zeta_l)\hat{B}^{\prime}_{j_l}(\zeta_l)\, d\zeta_l  \\
\end{array}
\end{equation}

For each integral in \eqref{eq:integrals} we define a quadrature
rule: we look for \begin{itemize}
\item points $\tilde{\xx}_{\boldsymbol{q}} = ( \tilde{x}_{1,q_1} ,
  \tilde{x}_{2,q_2}, \dots , \tilde{x}_{d,q_d} )$ with $q_l = 1,\dots
  n_{\mathsf{QP},l}$, with  $N_{\mathsf{QP}}$ is $ \# \left\{\tilde{\xx} \right\} =\prod_{l=1}^{d} n_{\mathsf{QP},l}$;
\item for each index $i_l=1,\dots, n_{\mathsf{DOF},l}; l=1,\dots, d$,  four quadrature rules such that:
  \begin{equation}
    \label{eq:quadrule}
\begin{gathered}
\mathbb{Q}^{(0,0)}_{i_l} (f) := \sum_{q_l = 1}^{n_{\mathsf{QP},l}} w^{(0,0)}_{l, i_l, q_l } f( \tilde{x}_{l,q_l} ) \approx \int_0^1 f(\zeta_l) \hat{B}_{i_l} (\zeta_l) d\zeta_l \,; \\
\mathbb{Q}^{(1,0)}_{i_l} (f) := \sum_{q_l = 1}^{n_{\mathsf{QP},l}} w^{(1,0)}_{l, i_l, q_l } f( \tilde{x}_{l,q_l} ) \approx \int_0^1 f(\zeta_l) \hat{B}_{i_l} (\zeta_l) d\zeta_l \,; \\
\mathbb{Q}^{(0,1)}_{i_l} (f) := \sum_{q_l = 1}^{n_{\mathsf{QP},l}} w^{(0,1)}_{l, i_l, q_l } f( \tilde{x}_{l,q_l} ) \approx \int_0^1 f(\zeta_l) \hat{B}^{\prime}_{i_l} (\zeta_l) d\zeta_l \,; \\
\mathbb{Q}^{(1,1)}_{i_l} (f) := \sum_{q_l = 1}^{n_{\mathsf{QP},l}} w^{(1,1)}_{l, i_l, q_l } f( \tilde{x}_{l,q_l} ) \approx \int_0^1 f(\zeta_l) \hat{B}^{\prime}_{i_l} (\zeta_l) d\zeta_l \, .
\end{gathered}
  \end{equation}
fulfilling  the exactness requirement:
\begin{equation} \label{eq:exact}
\displaystyle  \begin{array}{l}
\mathbb{Q}^{(0,0)}_{i_l} (\hat{B}_{j_l} ) =  \mathbb{I}^{(0,0)}_{l,i_l,j_l} \, \\
\mathbb{Q}^{(1,0)}_{i_l} (\hat{B}^{\prime}_{j_l} ) =  \mathbb{I}^{(1,0)}_{l,i_l,j_l} \, \\
\mathbb{Q}^{(0,1)}_{i_l} (\hat{B}_{j_l} ) =  \mathbb{I}^{(0,1)}_{l,i_l,j_l} \, \\
\mathbb{Q}^{(1,1)}_{i_l} (\hat{B}^{\prime}_{j_l} ) =  \mathbb{I}^{(1,1)}_{l,i_l,j_l} \,
\end{array} \,,\ \forall j_l \in \mathcal{I}_{l,i_l} \, .
\end{equation}

\end{itemize} 
For stability  we also require that the quadrature rules
$\mathbb{Q}^{(\cdot,\cdot)}_{i_l} $ have support included in the
support of $\hat{B}_{i_l}$, that is  
\begin{equation} \label{eq:support}
q_l  \notin \mathcal{Q}_{l,i_l}  \Rightarrow w^{(\cdot,\cdot)}_{l, i_l , q_l} =0  \ .
\end{equation}
where $\mathcal{Q}_{l,i_l}:= \left\{ q_l\in {1,\dots,
    n_{\mathsf{QP},l}} \text{ s.t. } \tilde{x}_{l,q_l} \in
  \mbox{supp}\left(\hat{B}_{i_l}\right)
\right\} $; recall that here the support of a function is considered an open set. Correspondingly, we introduce  the set of multi-indexes $\mathcal{Q}_{\ii }
:= \prod_{l=1}^d \mathcal{Q}_{l,i_l} $.

Once the points $\tilde{\xx}_{\boldsymbol{q}} $ are
fixed, the quadrature rules have to be determined by the exactness
requirements, that are a system of  linear equations of the unknown
weights (each of the \eqref{eq:exact}).
For that  we require
\begin{equation}\label{eq:numb_nodes}
\# \mathcal{Q}_{l,i_l} \ge \# \mathcal{I}_{l,i_l} \ .
\end{equation}
See Remark \ref{rem:well-posedness-for-weights} for a discussion on the
well-posedness of the linear systems for the weights.

\begin{figure}
\begin{center}
 \includegraphics[width=0.45\textwidth] {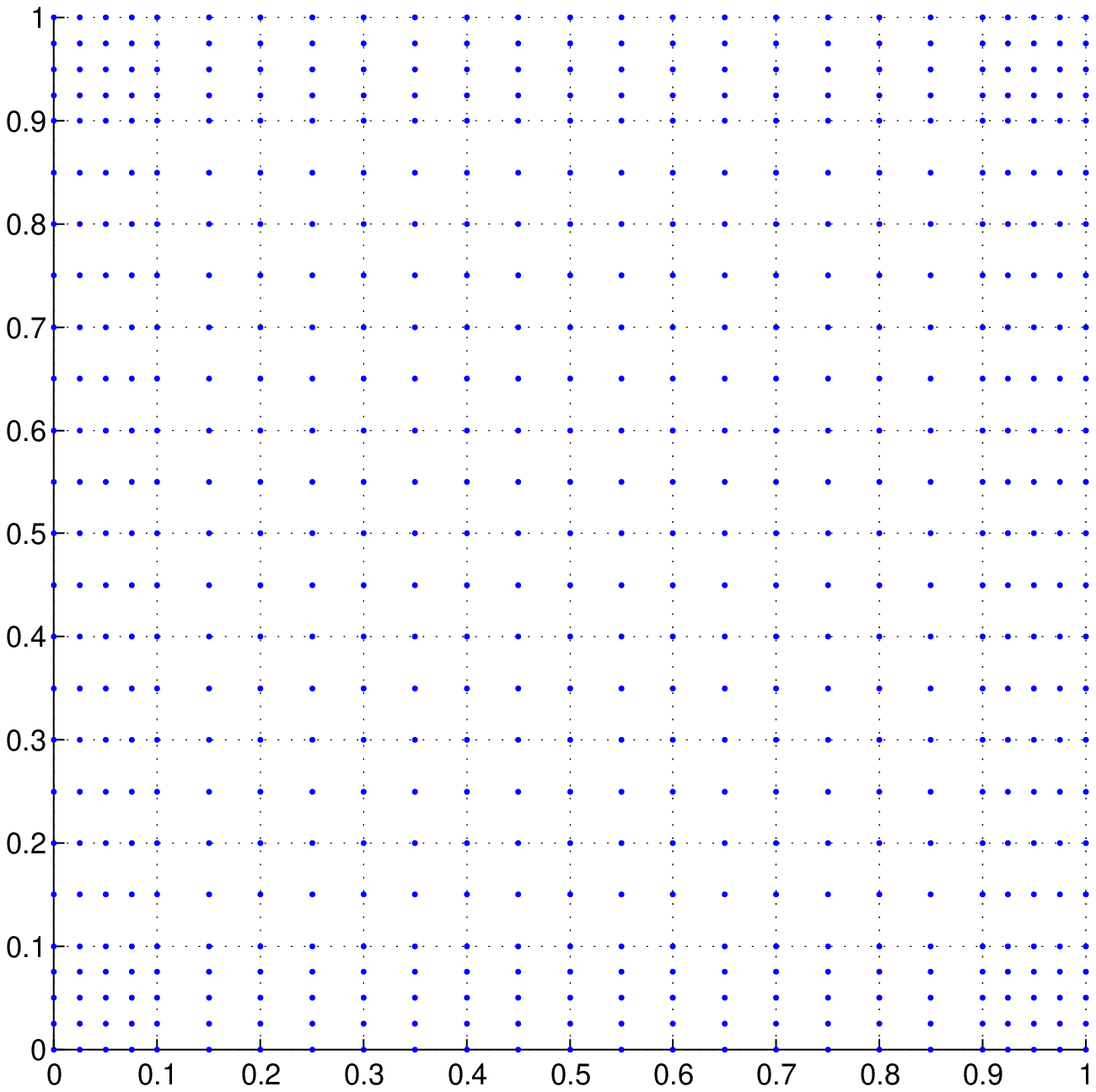} 
 \includegraphics[width=0.45\textwidth] {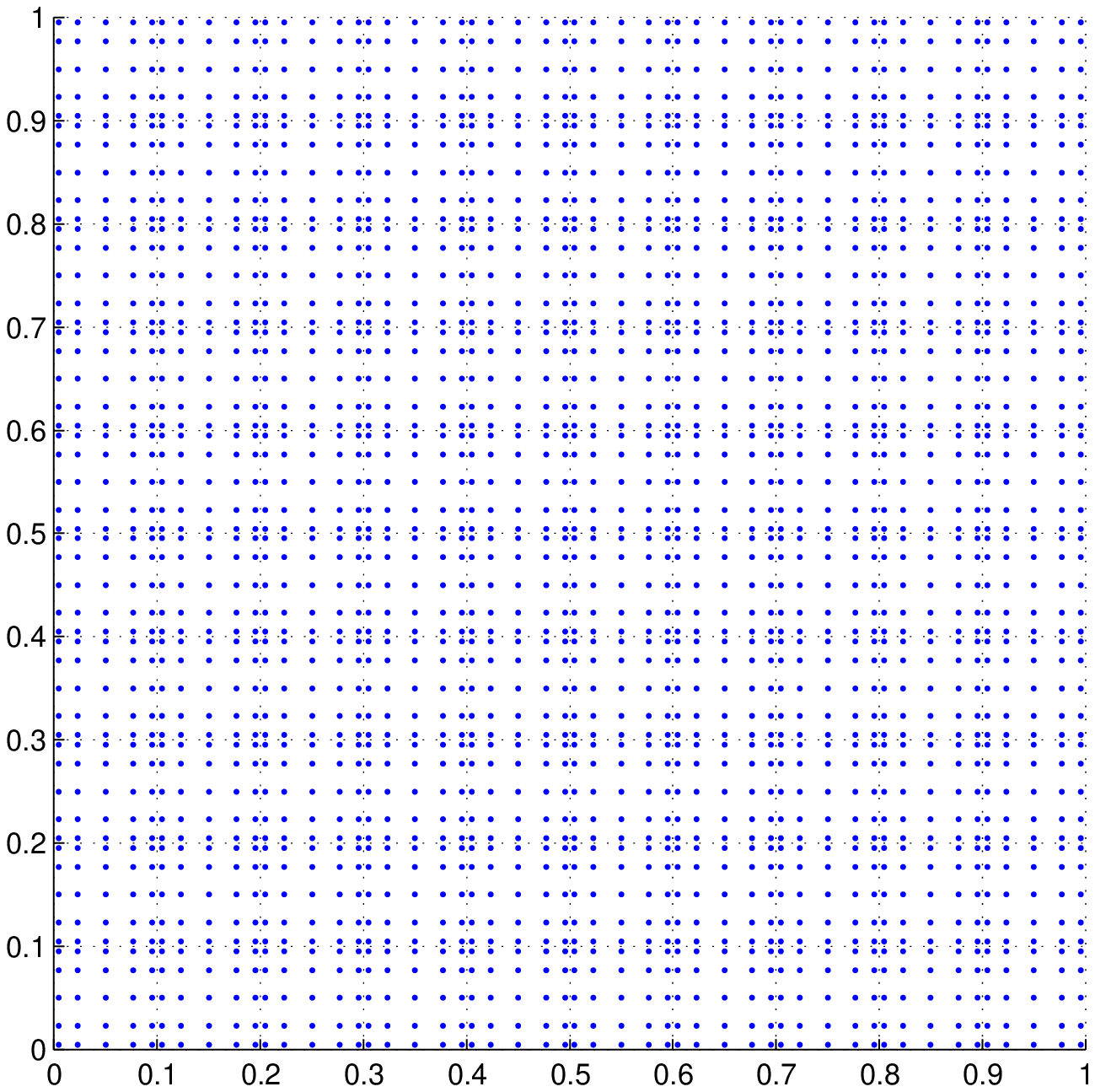} \\
  \includegraphics[width=0.45\textwidth] {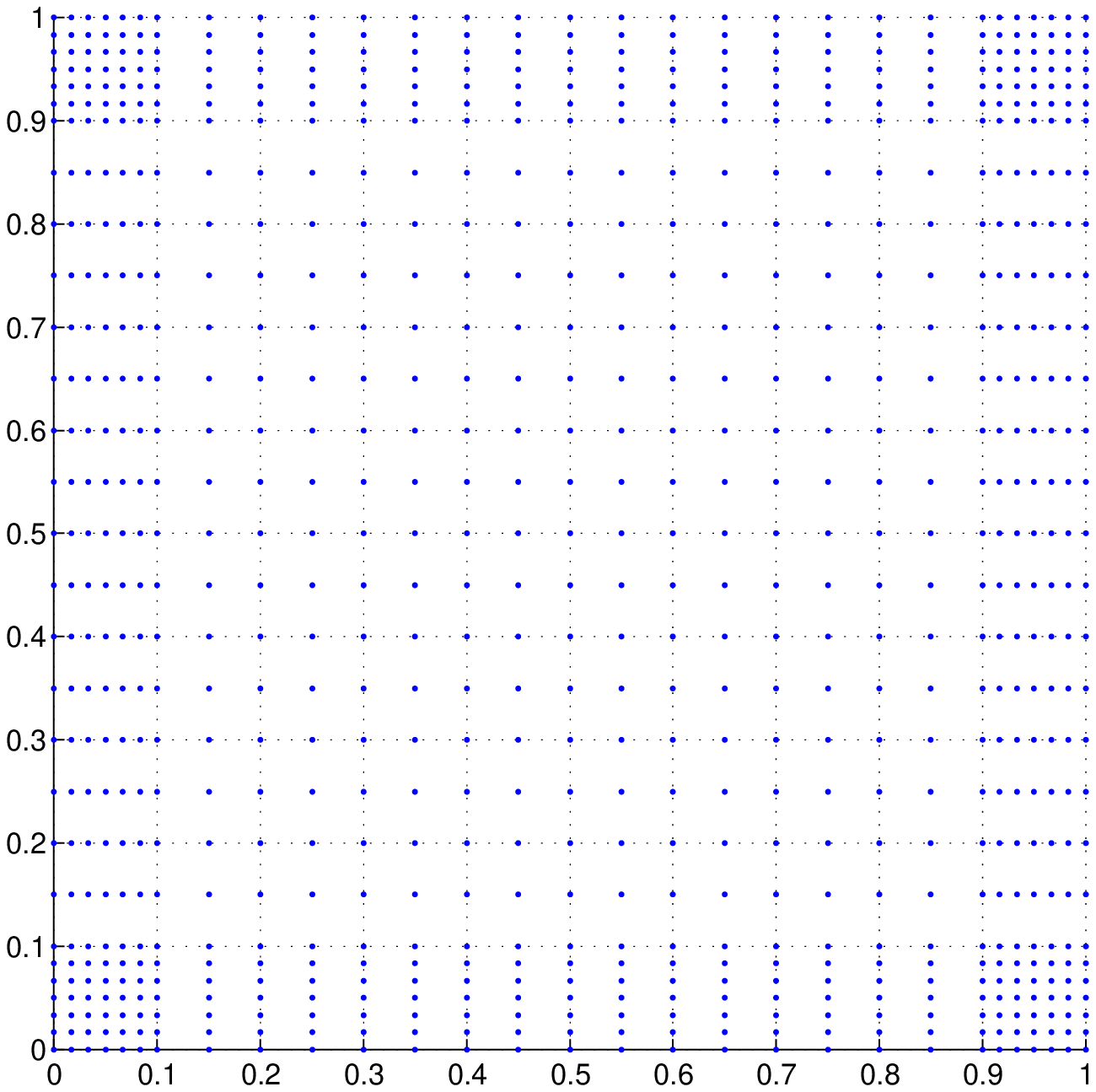} 
 \includegraphics[width=0.45\textwidth] {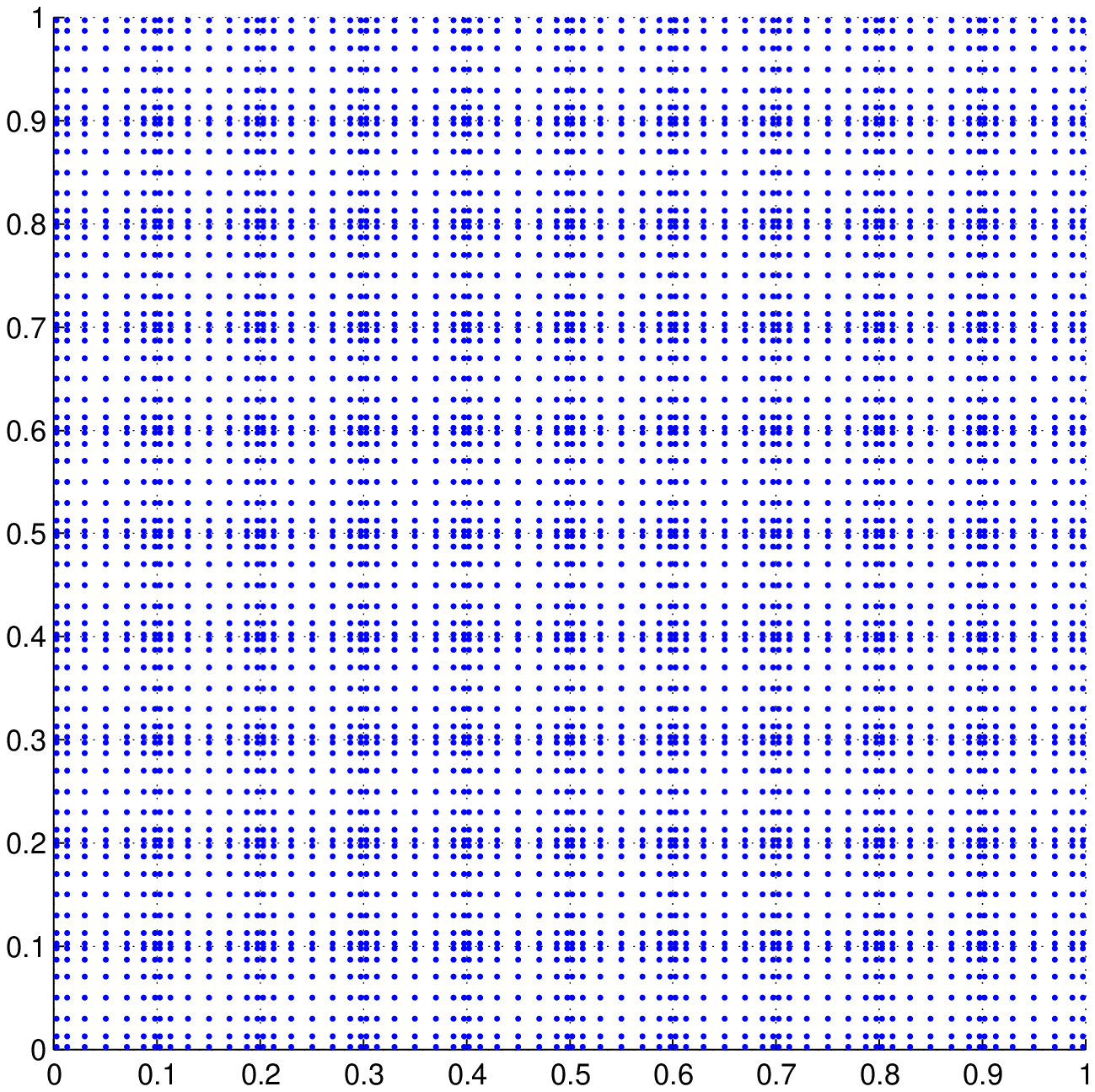}
\end{center}
\caption{Comparison between quadrature rules. On the first line, the
  quadrature  points needed for the case $p=4$, $n_{\mathsf{EL}}=10$,
  $d=2$, in the second line the case $p=6$. On the left panel the
  proposed WQ rule, on the right the SGQ  rule. } \label{fig:nodes}
\end{figure}

\begin{remark}[The choice of quadrature points]\label{Rem:points}
The construction of a global grid of quadrature points is done in
order to save computations. For the case of maximum $C^{p-1}$
regularity considered here, our
choice for quadrature points  is endpoints (knots) and midpoints of
all internal knot-spans, while for  the boundary knot-spans
(i.e. those that are adjacent 
to the boundary of the parameter domain $\hat \Omega$) we take $p+1$ equally spaced 
points. Globally    $N_{\mathsf{QP}} \approx 2^d N_{\mathsf{EL}}
= O(N_{\mathsf{DOF}} ) $ considering only the dominant term (remember
that $n_{\mathsf{EL},l} \gg p$).  In
Figure \ref{fig:nodes} we plot the quadrature points grid, and a
comparison is made with respect to element-by-element standard Gaussian
quadrature (SGQ) points. 
\end{remark}
\begin{remark}[Computation of quadrature weights]\label{rem:well-posedness-for-weights}
Given the quadrature points, the quadrature weights  are selected in
order to fulfil \eqref{eq:exact}--\eqref{eq:support}. When $ \# \mathcal{Q}_{l,i_l} >\# \mathcal{I}_{l,i_l} 
 $ the  quadrature weights are not uniquely given from
 \eqref{eq:exact}--\eqref{eq:support}  and are selected by a minimum
 norm condition. In all  cases  with our choice of quadrature points and thanks to the
Schoenberg-Whitney interpolation theorem \cite[Theorem 2 \P
XIII]{DeBoor:1978} we can guarantee that the quadrature weights
fulfilling the above conditions exist. 

\end{remark}
\begin{remark}[Alternative choices of the quadrature points]\label{rem:alt-quad-points}
% The proposed method distinguishes the quadrature to be used along each line of the matrices involved in the discrete counterpart of the PDE, thus a selection has to be made when applying quadrature. For this reason, we have made the choice to present four different quadrature rules for the four possible integrals arising, but one is free to choose instead two quadrature rule that are exact both for B-splines and their derivatives.
If there is no need for a global grid of quadrature points
(e.g., the cost of calculation of the coefficients is negligible), it
is possible to have quadrature points that  depend on the index $i_l$,
as for the weights. Then,  one can construct weighted Gaussian
quadrature rules (see e.g. \cite{Calabro:2015}) in order to
minimize the number of quadrature points associated to each row of the
stiffness matrix.
\end{remark}

\section{Pseudo-codes and computational cost}\label{Sect:Pseudocodes}

In order to simplify the FLOPs count,  we assume
$n_{\mathsf{EL}}:=n_{\mathsf{EL},1}= n_{\mathsf{EL},2}= \dots
=n_{\mathsf{EL},d} $ and
$n_{\mathsf{DOF}}:=n_{\mathsf{DOF},1}=n_{\mathsf{DOF},2}=\dots =
n_{\mathsf{DOF},d}$. We then have $N_{\mathsf{EL}} = n_{\mathsf{EL}}^d
$, $N_{\mathsf{DOF}} = n_{\mathsf{DOF}}^d $, etc.
We consider the case of maximum regularity $r=p-1$ and
$n_{\mathsf{EL}} \gg p$ that implies $n_{\mathsf{EL}} \approx
n_{\mathsf{DOF}}$. We recall that:
\begin{equation}\label{eq:prelim_costs}
\begin{array}{r}
\# \mathcal{K}_{l,i_l} \leq p+1 \\\# \mathcal{I}_{l,i_l} \leq 2p+1 
\end{array} 
\end{equation}
With our choice for the quadrature points, the previous two imply $\# \mathcal{Q}_{l,i_l}
\le 2p+1$.

We first collect all the initialisations needed in Algorithm
\ref{algo_iniz}. It is not necessary to  precompute these quantities 
-and in most architectures access to stored  data is costly- but this
used here for FLOPs  evaluation.

\begin{algorithm} \LinesNumbered \SetKwInOut{Input}{Input}
\Input{Quadrature points $\tilde{\xx}_{\boldsymbol{q}}$ as in Remark \ref{Rem:points} }
\For{$l=1,\dots, d$}{
\For{ $i_l=1,\dots, n_{\mathsf{DOF},l} $}{ 
Evaluate $\mathbb{B}_{l,i_l,q}^{(0)} := B_{i_{l}}( \tilde{x}_{l,q}) \, \forall q \in \mathcal{Q}_{l,i_l} $, store  $\mathbb{B}_{l,i_l}^{(0)}\in \R^{ n_{\mathsf{QP},l} } $ \;
Evaluate $\mathbb{B}_{l,i_l,q}^{(1)} := B^{\prime}_{i_{l}}( \tilde{x}_{l,q}) \, \forall q \in \mathcal{Q}_{l,i_l}$, store  $\mathbb{B}_{l,i_l}^{(1)}\in \R^{ n_{\mathsf{QP},l} } $ \;
}
\For{ $i_l=1,\dots, n_{\mathsf{DOF},l} $}{ 
\For{ $j_l \in \mathcal{I}_{l,i_l} $}{ 
Calculate $\mathbb{I}^{(0,0)}_{l,i_l,j_l}$, $\mathbb{I}^{(1,0)}_{l,i_l,j_l}$, $\mathbb{I}^{(0,1)}_{l,i_l,j_l}$, $\mathbb{I}^{(1,1)}_{l,i_l,j_l}$  as defined in \eqref{eq:integrals} \;
}
}

\For{$m=1,\dots, d$}{
Evaluate $c_{l,m}(\zeta) $ of equation \eqref{eq:coef_stiff} on points $\tilde{\xx}_{\boldsymbol{q}}$ \;
}
}
Evaluate $c(\zeta)$ on points $\tilde{\xx}_{\boldsymbol{q}}$\;
\caption{Initializations}\label{algo_iniz}
\end{algorithm}

  Then we can count operations in Algorithm \ref{algo_iniz}:
\begin{enumerate}
\item[(i)] Evaluations  of B-splines  reported on lines 3--4 can be
  done in  $\frac{1}{2} p^2$ FLOPs each. They are repeated  $d\,  n_{\mathsf{DOF}} \,  \# \mathcal{Q}_{l,i_l}  $ times so that this part costs $O\left( p^3 \,  n_{\mathsf{DOF}} \right)$ FLOPs.
\item[(ii)] The calculation of integrals \eqref{eq:integrals} on line
  8 needs to be done in an exact manner. The calculation of the exact
  integral of products of B-splines has a vast literature
  \cite{vermeulen1992integrating}, however, closed forms are available only
  in some particular cases. For this reason we consider here  the
  usual element-wise Gaussian quadrature.  
The evaluation of B-splines and their derivatives cost, as reported before, $\frac{1}{2} p^2$ FLOPs for each point. Counting all Gaussian point, the cost is $\approx d\,  p^2$ evaluations of each of the $d\,  n_{\mathsf{DOF}}$ univariate basis functions, thus costs $O(p^4 \,  n_{\mathsf{DOF}})$ FLOPs.\\ 
The computation of each of the integrals has the cost of a summation on $\approx p^2$ terms; and the four calculations are done $ d\,  n_{\mathsf{DOF}}  \,  \# \mathcal{I}_{l,i_l} $ times so that this costs  $O\left(p^3 \,  n_{\mathsf{DOF}}\right)$ FLOPs.
\item[(iii)] The evaluations of the $ (d^2+1)$ functions
  $c_{l,m} $ and $c$
  on lines 12 and 15 have to be performed at the
  $N_{\mathsf{QP}} = (n_{\mathsf{QP}})^d$ quadrature points. The
  actual cost depends on the evaluation cost of $c_{l,m} $ and
  $c$. If these coefficients are  obtained by  $O(p^d)$ linear
  combinations of B-spline values (or derivatives), and each multivariate B-spline value is computed
  from  multiplications of  univariate B-spline values, the total cost
   is  $C(d)p^d$ per quadrature point and in total  $O\left( p^d \, N_{\mathsf{QP}} \right)$ FLOPs. 
\end{enumerate} 
The leading cost of Algorithm 1 for $d\geq 2$  is $O\left( p^d \,  N_{\mathsf{QP}} \right)$. 

\begin{algorithm} \LinesNumbered \SetKwInOut{Input}{Input}
\Input{Quadrature points $\tilde{\xx}_{\boldsymbol{q}}$, B-spline evaluations $\mathbb{B}_{l,i_l,\boldsymbol{q}}^{(\cdot )}$, Integrals $\mathbb{I}^{(\cdot ,\cdot )}_{l,i_l,j_l}$}
\For{$l=1,\dots, d$}{
\For{ $i_l=1,\dots, n_{l,\mathsf{DOF}}$}{
Calculate $w^{(0,0)}_{l, i_l, \mathcal{Q}_{l,i_l}}$ as  (minimum
Euclidean norm) solution of $\mathbb{B}_{l,\mathcal{I}_{l,i_l},\mathcal{Q}_{l,i_l}}^{(0)}  w^{(0,0)}_{l, i_l, \mathcal{Q}_{l,i_l} } = \mathbb{I}^{(0 ,0 )}_{l,i_l,\mathcal{I}_{l,i_l}} $ \;
Calculate $w^{(1,0)}_{l, i_l, \mathcal{Q}_{l,i_l}}$ as  (minimum
Euclidean norm)  solution of $\mathbb{B}_{l,\mathcal{I}_{l,i_l}, \mathcal{Q}_{l,i_l}}^{(1)}  w^{(1,0)}_{l, i_l, \mathcal{Q}_{l,i_l} } = \mathbb{I}^{(1 ,0 )}_{l,i_l,\mathcal{I}_{l,i_l}} $ \;
Calculate $w^{(0,1)}_{l, i_l, \mathcal{Q}_{l,i_l}}$ as (minimum
Euclidean norm)  solution of $\mathbb{B}_{l,\mathcal{I}_{l,i_l},\mathcal{Q}_{l,i_l} }^{(0)}  w^{(0,1)}_{l, i_l, \mathcal{Q}_{l,i_l} } = \mathbb{I}^{(0 ,1 )}_{l,i_l,\mathcal{I}_{l,i_l}} $ \;
Calculate $w^{(1,1)}_{l, i_l, \mathcal{Q}_{l,i_l}}$ as  (minimum
Euclidean norm) solution of $\mathbb{B}_{l,\mathcal{I}_{l,i_l}, \mathcal{Q}_{l,i_l}}^{(1)} w^{(1,1)}_{l, i_l, \mathcal{Q}_{l,i_l} } = \mathbb{I}^{(1 ,1 )}_{l,i_l,\mathcal{I}_{l,i_l}} $ \;
 } 
 }
\caption{Construction of univariate WQ rules}\label{algo_quad}
\end{algorithm}

In Algorithm \ref{algo_quad} we summarize the operations needed for
the construction of the univariate  WQ rules. Each calculation in lines
3--6 consists in the resolution of a linear system of dimension
$\approx (2p)^2$ that is possibly   under-determined.  The cost of
these computations in any case negligible since it is proportional to $n_{\mathsf{DOF}} $.
% Notice, in particular, that the matrix inversions are of lower degree, being of order $O \left( (3p)^3\,  \sum_{l=1}^d    n^{(l)}_{\mathsf{DOF}}   \right) $.  
%\begin{remark}[Connection matrix]
%\end{remark}

% \begin{remark}[Regular meshes]
% In our construction, we allow the construction of tensor product spline spaces with different knot vectors along each direction and general non-uniform spacing. Nevertheless we explicitly notice that in these remarkable cases we can explore symmetries in order to avoid redounding calculations. In particular:
% \begin{itemize}
% \item if $\Xi_{l_1}= \Xi_{l_2} \forall l_1, l_2 \in \{1,\ldots, d\} $ then the calculation of the quadrature rules has to be done only one time and used $d$ times;
% \item if $\xi_{j+1}-\xi_j =\lambda \forall j=0,\ldots, n_{\mathsf{EL}}-1$  (uniform patch) then most of the basis functions are translation copies of the cardinal B-spline thus the quadrature rule with respect to this function can be made only one time;
% \item the case of patches with geometrically-scaled lengths can be explicitly taken in account as done in \cite{Auricchio:2012}.
% \end{itemize}
% \end{remark}

\section{Formation of the mass matrix}
When all the quadrature rules are available we can write the
computation of the approximate mass matrix following
\eqref{eq:Mass_approx}. 
Similar formulae and algorithms can be written for the stiffness
matrix starting from equation \eqref{eq:stiffness-3d}.

The mass matrix formation algorithm is mainly a loop over all rows
$\ii$,  for each $\ii$ we consider the calculation of  
\begin{equation}\label{eq:WQ-quadrature}
\tilde{m}_{\ii,\jj}= \mathbb{Q}^{(0,0)}_{\ii} \left(
  \hat{B}_{\jj}(\zzeta) c(\zzeta)  \right), \quad \forall \jj\in \mathcal{I}_{\ii}.
\end{equation}
% \begin{algorithm} \LinesNumbered \SetKwInOut{Input}{Input}
% \Input{Quadrature rules, evaluations of coefficients }
% \For{$\ii=1,\dots, N_{\mathsf{DOF}}$}{
% Load (or construct) quadrature rules $\mathbb{Q}^{(0,0)}_{i_l}\,,\ l=1,\dots, d$ \;
% \For{$\jj\in \mathcal{I}_{\ii} $}{
% Compute $\tilde{m}_{\ii,\jj}= \mathbb{Q}^{(0,0)}_{\ii} \left( \hat{B}_{\jj}(\zzeta) c(\zzeta)  \right) $\;
% }
% }
% \caption{Construction of Mass matrix}\label{algo_mass_v1}
% \end{algorithm}
The computational cost of \eqref{eq:WQ-quadrature} is minimised by a sum
factorization approach, which is explained below. 
%In Algorithm \ref{algo_mass} we report the sum-factorization version of Algorithm  \ref{algo_mass_v1}. 

If we substitute \eqref{eq:quadrule} into \eqref{eq:Mass_approx} we obtain the following sequence of nested summations:

\begin{equation} \label{eq:massquadrature} \tilde m_{\ii, \jj} = \sum_{q_1 \in \mathcal{Q}_{1,i_1}} w^{(0,0)}_{1,i_1,q_1} \hat B_{j_1}(\tilde x_{1,q_1}) 
\left( \sum_{q_2 \in \mathcal{Q}_{2,i_2}} 
%w^{(0,0)}_{2,i_2,q_2} B_{j_2}(\tilde x_{2,q_2}) \left( 
\ldots 
%\left(
\sum_{q_d \in \mathcal{Q}_{d,i_d}} w^{(0,0)}_{d,i_d,q_d} \hat B_{j_d}(\tilde x_{d,q_d}) c\left( x_{1,q_1}, \ldots, x_{d,q_d} \right) 
%\right) 
\right) \end{equation}

To write \eqref{eq:massquadrature} in a more compact form, we
introduce the notion of matrix-tensor product \cite{Kolda:2009}. Let
$\mathcal{X} = \left\{x_{k_1, \ldots, k_d}\right\} \in
\mathbb{R}^{n_{1} \times  \ldots \times  n_{d} }$ be a $d-$dimensional
tensor, and let $m \in \left\{1, \ldots, d\right\}$. The $m-$mode
product of $\mathcal{X}$ with a matrix $A = \left\{a_{i,j}\right\} \in
\mathbb{R}^{ t \times n_m}$, denoted with $\mathcal{X} \times _{m} A $,
is a tensor of dimension $n_1 \times  \ldots \times  n_{m-1} \times  t
\times  n_{m+1} \times  \ldots \times  n_d$, with  components 
$$ \left(\mathcal{X} \times_m A \right)_{k_1, \ldots, k_d} = \sum_{j =
  1}^{ n_m} a_{k_m, j} \; x_{k_1, \ldots k_{m-1}, j, k_{m+1}, \ldots k_d}. $$
We emphasize that such computation requires $2 t \prod_{l = 1}^d n_l$ FLOPs.

For $l= 1, \ldots, d$ and $i_l = 1, \ldots, n_{\mathsf{DOF},l}$ we define the matrices
$$ \mathrm{B}^{(l,i_l)} = \left( \hat{B}_{j_l} (x_{l,q_l}) \right)_{j_l \in \mathcal{I}_{l,i_l}, q_l \in \mathcal{Q}_{l,i_l}}, \qquad \mathrm{W}^{(l,i_l)} = \mbox{diag}\left(\left( w^{(0,0)}_{l,i_l,q_l} \right)_{q_l \in \mathcal{Q}_{l,i_l}} \right), $$
where $\mbox{diag}(v)$ denotes the diagonal matrix obtained by the vector $v$. We also define, for each index $\ii$, the $d-$dimensional tensor
$$ \mathcal{C}_{\ii} = c( \tilde{\xx}_{\mathcal{Q}_{\ii} } ) = \left(c(\tilde{x}_{1,q_1}, \ldots, \tilde{x}_{d,q_d})\right)_{q_1 \in \mathcal{Q}_{1,i_1}, \ldots, q_d \in \mathcal{Q}_{d,i_d}} $$
%
%We now fix $\ii = 1, \ldots,$ and consider \eqref{eq:massquadrature} for every $\jj \in \mathcal{I}_{\ii}$.
 Using the above notations, we have

\begin{equation} \label{eq:sumfactorization} \tilde m_{\ii, \mathcal{I}_{\ii}} = \mathcal{C}_{\ii} \times_d \left(\mathrm{B}^{(d,i_d)} \mathrm{W}^{(d,i_d)}\right) \times_{d-1} \ldots \times_1 \left(\mathrm{B}^{(1,i_1)} \mathrm{W}^{(1,i_1)}\right). \end{equation}

%Note that, when $d=2$, we have $\cal{X} \times_1 A = \cal{X} A^T $ and $\cal{X} \times_2 A = A \cal{X} $.

Since with our choice of the quadrature points $\#
\mathcal{Q}_{l,i_l}$ and $\# \mathcal{I}_{l,i_l}$ are both $O(p)$, the
computational cost associated with \eqref{eq:sumfactorization} is
$O(p^{d+1})$ FLOPs. Note that $\tilde m_{\ii, \mathcal{I}_{\ii}}$
includes all the nonzeros entries of the $\ii-$th row of
$\tilde{\mathbb{M}}$. Hence if we compute it for each $\ii =
1,\ldots,N_{\mathsf{DOF}}$ the total cost amounts to
$O(N_{\mathsf{DOF}} \,  p^{d+1})$ FLOPs. This approach is summarized in Algorithm \ref{algo_mass}.

We remark that writing the sums in \eqref{eq:massquadrature} in terms of matrix-tensor products as in \eqref{eq:sumfactorization} is very useful from an implementation viewpoint.
Indeed, in interpreted languages like MATLAB (which is the one used in the experiments of the next section), 
it is crucial to avoid loops and vectorize (in our case, tensorize) the operations, in order to obtain an efficient implementation of an algorithm; see also the discussion in \cite{Cuvelier:2015}.
In particular, each matrix-tensor product in
\eqref{eq:sumfactorization} is computed via a simple matrix-matrix
product, which is a BLAS level 3 operation and typically yields high efficiency on modern computers.

%By a sum-factorization-type recursion, we can obtain an efficient procedure for the previous summations, see Algorithm \ref{algo_mass}. In this case we work with a fixed index $\ii$ at the time and we construct the whole row of the matrix. 

%In order to calculate the computational cost, we note that for each $\ii$ the calculation at line 9 costs $2\times \# \mathcal{Q}_{l,i_l}$ FLOPs. Moreover this calculation is performed more times when the index $l$ increases, and when $l=d$ is done $ \# \mathcal{I}_{\ii} \approx (2p)^d $ times. Then, the computational cost can be estimated in $O\left(2^{d+2} \times N_{\mathsf{DOF}}\times p^{d+1} \right)$ FLOPs. 

\begin{algorithm} \LinesNumbered \SetKwInOut{Input}{Input}
\Input{Quadrature rules, evaluations of coefficients }
\For{$\ii=1,\dots, N_{\mathsf{DOF}}$}{
Set $\mathcal{C}^{(0)}_{\ii} := c( \tilde{\xx}_{\mathcal{Q}_{\ii} } ) $\;
\For{ $l=d,d-1,\ldots , 1$}{
Load the quadrature rule $\mathbb{Q}^{(0,0)}_{i_l}$ and form the matrices $\mathrm{B}^{(l,i_l)}$ and $\mathrm{W}^{(l,i_l)}$\; 
Compute $\mathcal{C}^{(d+1-l)}_{\ii} = \mathcal{C}^{(d-l)}_{\ii} \times_l \left(\mathrm{B}^{(l,i_l)} \mathrm{W}^{(l,i_l)}\right) $\;
}
Store $\tilde{m}_{\ii, \mathcal{I}_{\ii} } = \mathcal{C}^{(d)}_{\ii} $\;
}
\caption{Construction of mass matrix by sum-factorization}\label{algo_mass}
\end{algorithm}

\section{Numerical tests}
In this section, in order to evaluate numerically the behavior of the proposed procedure we present some numerical tests. 
First, in Section \ref{Sect5.1} we consider the solution of a 1D
problem where we see that the  application of our  row-loop WQ-based
algorithm  leads to optimal order of convergence.  
Then, in Section \ref{Sect5.2} we measure the performance of the
algorithm. We consider there  the formation
of mass matrices in 3D. 
The results for all cases refers to a 
Linux workstation equipped with Intel i7-5820K processors running at 3.30GHz, and with 64 GB of RAM.
% Linux workstation equipped with Xeon E5-2470 processors (running at 2.3 GHz), \B  
 The row-loop WQ-based algorithm  is potentially  better suited for a
parallel implementation than the standart element-wise  SGQ-based
algorithm, however we benchmark here sequential execution and use only one core for the simulations.

\subsection{Convergence of approximate solution in 1D} \label{Sect5.1}
\begin{figure}
\begin{center}
 \includegraphics[width=0.49\textwidth] {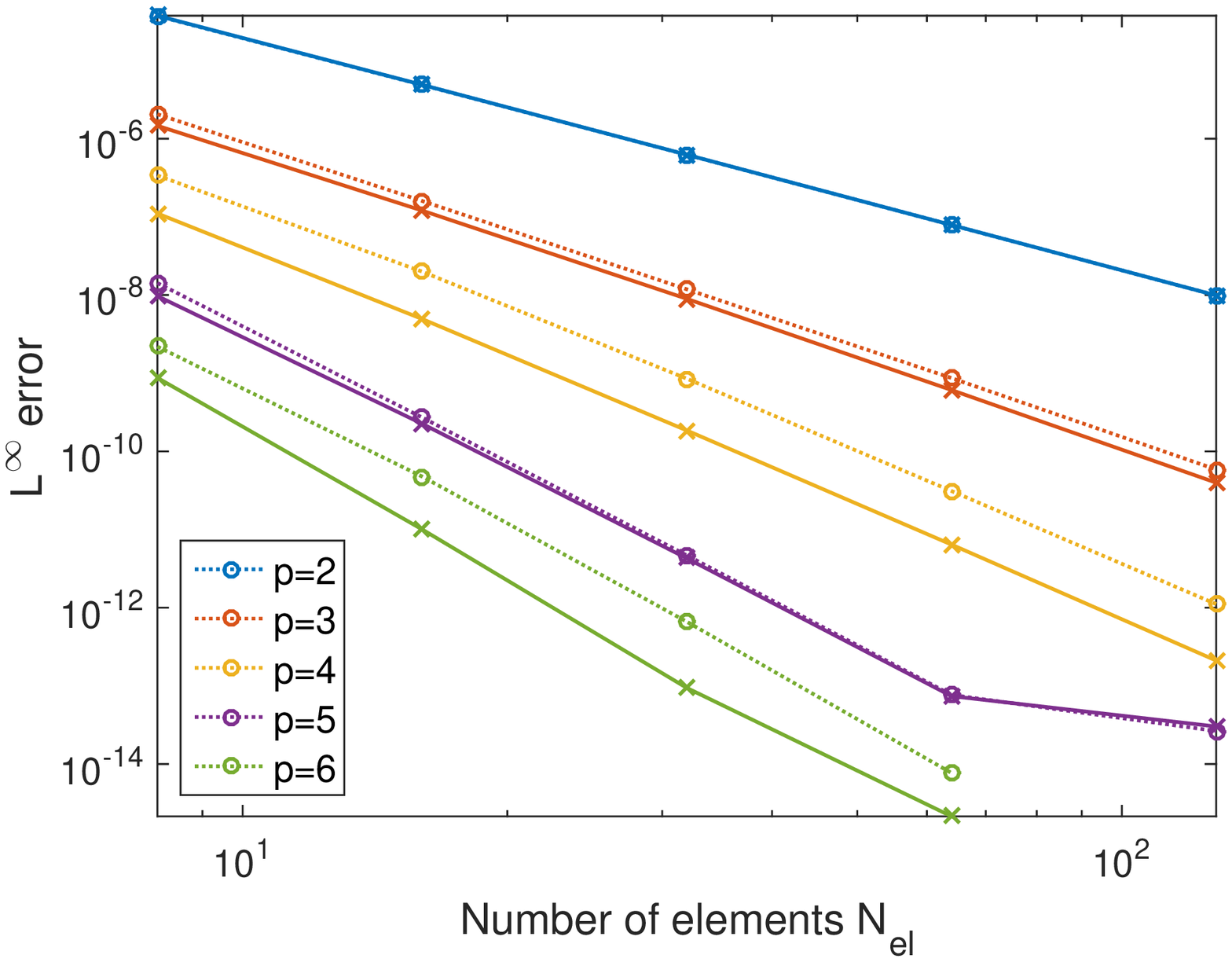} 
 \includegraphics[width=0.49\textwidth] {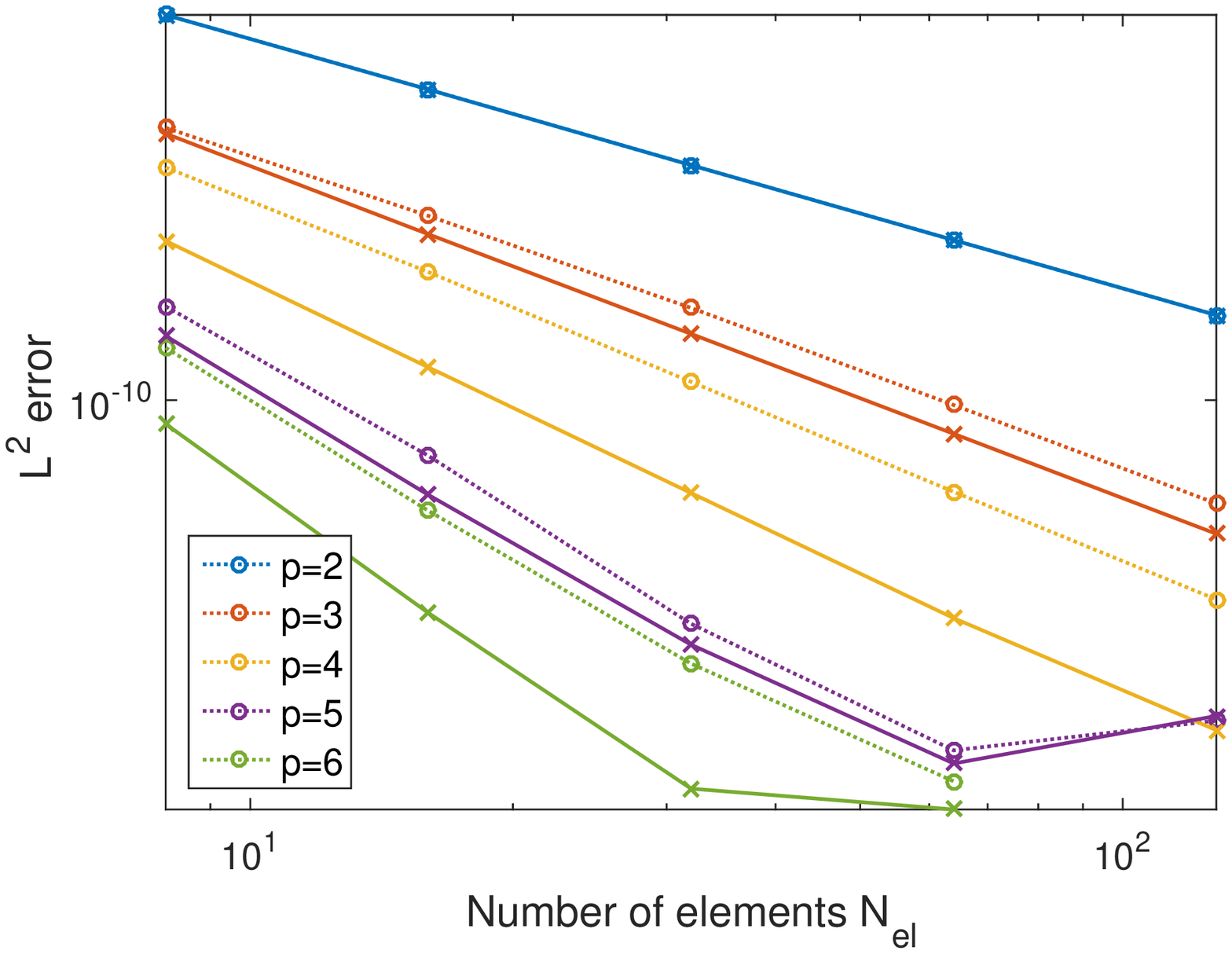} \\ \includegraphics[width=0.49\textwidth] {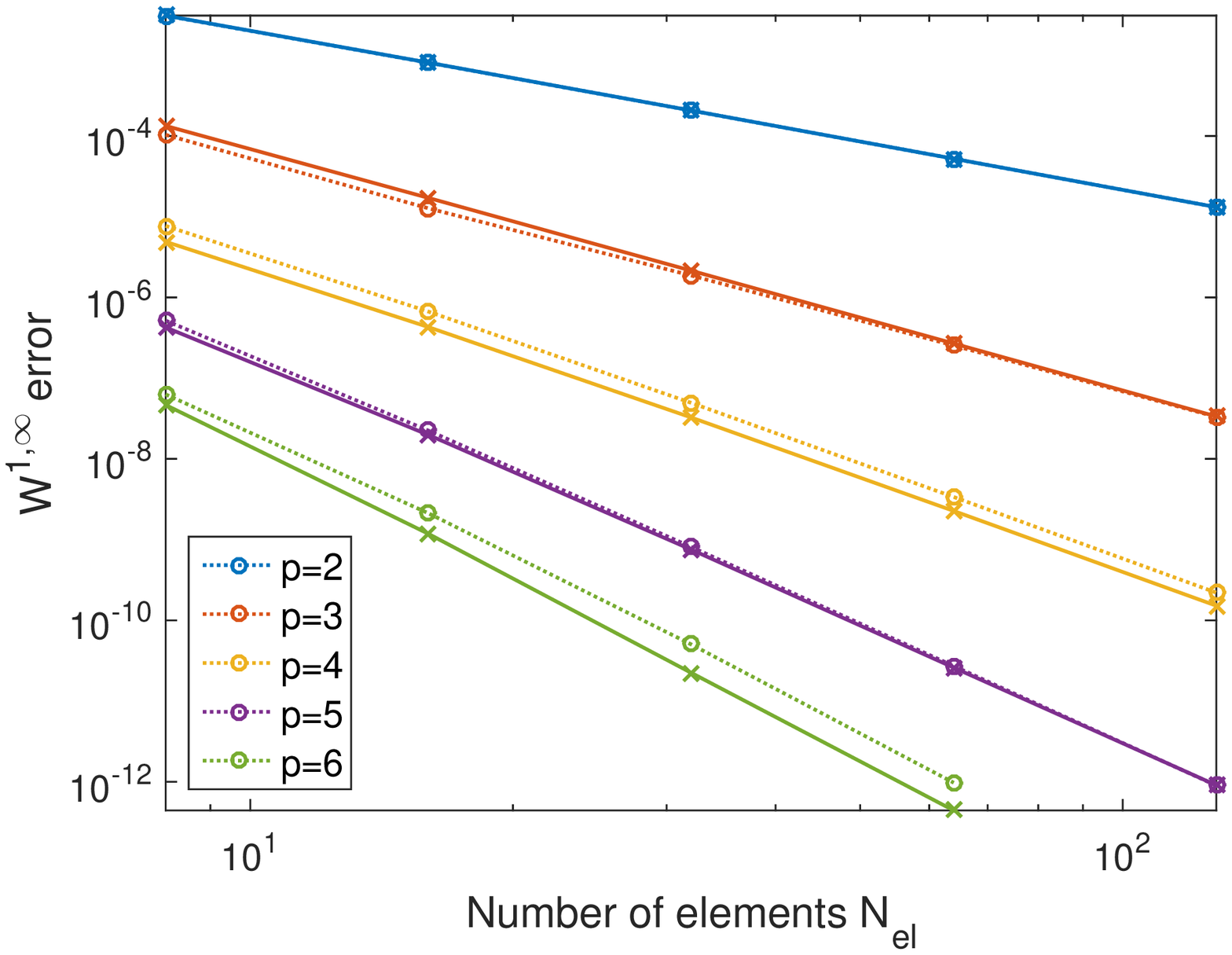} \\
\end{center}
\caption{Convergence history plot. We report errors in  $L^{\infty}, L^2$ and $H^1$
  norm for the solution of problem \eqref{poisson1D} by
  Galerkin based isogeometric analysis with WQ (Algorithm
  \ref{algo_iniz}--\ref{algo_mass}) for various degrees $p$ in dotted
  lines. As reference, the solid lines  refer to the same calculation
  made with element-wise SGQ.  Optimal convergence rate
  is achieved in all cases. SGQ is sligtly  more accurate  for even
  degrees $>2$ in  $L^2$ and $L^{\infty}$ norms. \B } \label{figure:1}
\end{figure}

As a test with known solution we consider the following:
\begin{equation}\label{poisson1D}
\left\{
\begin{aligned}
u^{\prime \prime} + u(x) = \frac{5 \mbox{exp}(2x)-1}{4} &\quad \text{on} \qquad [0,\pi/6] \\
u(0) = 0, u(\pi/6)= \mbox{exp}( \pi/3 )/2 & 
\end{aligned}
\right.
\end{equation}
We compare the numerical solution in the parametric domain, using the geometric transformation $t= 2 sin(x)$, with the exact one $u(x)=\frac{\mbox{exp}(2x)-1}{4}$. Then we calculate point-wise absolute error, integral error and energy error  - namely $L^{\infty}, L^2$ and $H^1$ norms - with varying spline degree $p$. Figure \ref{figure:1} illustrates that the construction of the matrices with the proposed procedure does not effect the overall convergence properties of the Galerkin method, as it can be seen by comparing the convergence curves with those obtained using Gaussian quadrature. 

\begin{remark}[Convergence rates]\label{rem:symmetry} 
As already noted in Section \ref{Sect:Weighted}, WQ  does not preserve
symmetry, that is in general $\tilde m_{\ii, \jj} \neq \tilde m_{\jj, \ii} $ even if $
m_{\ii, \jj} =  m_{\jj, \ii} $. The lack of symmetry did not cause any
deterioration of the order of convergence in energy and lower-order
norms in our numerical benchmarking, see Figure \ref{figure:1}.  This
is an  important and interesting behaviour that deserves
further study. We remark that  the lack of symmetry occurs also for  collocation isogeometric
schemes \cite{Schillinger:2013}, where however convergence rates are
suboptimal. \B
\end{remark}

\subsection{Time for the formation of matrices} \label{Sect5.2}
\begin{figure}
\begin{center}
 \includegraphics[width=0.6\textwidth]{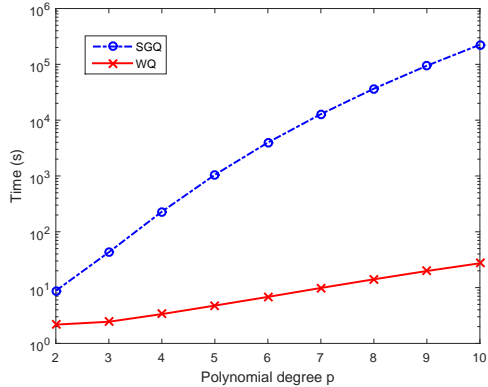} %Time_nel
\end{center}
\caption{ Time for mass matrix assembly in the
  framework of isogeometric-Galerkin method with maximal regularity on
  a single patch domain of $20^3$ elements. The comparison is between the
  WQ approach proposed (Algorithm \ref{algo_iniz}--\ref{algo_mass})
  and SGQ as implemented in GeoPDEs 3.0. } \label{figure:2}
\end{figure}

In this section we report CPU time results for the
formation on a single patch domain of mass matrices. Comparison is
made with  GeoPDEs, the optimized but SGQ-based MATLAB
isogeometric library developed by Rafael V\'azquez,  see
\cite{de2011geopdes,geopdesv3}. In Figure \ref{figure:2} we plot the
time needed for the mass matrix formation  up to degree $p=10$ with
$N_{\mathsf{DOF}}= 20^3$. The tests confirm the superior
performance of the proposed  row-loop WQ-based algorithm vs SGQ. In the case $p=10$ GeoPDEs takes more than 62 hours to form the mass matrix while the proposed algorithm needs only 27 sec, so that the use high degrees is possible with WQ.

\begin{remark}[Sparse implementation of WQ]
Clearly we exploit sparsity in our MATLAB implementation: we compute
all the nonzero entries of $\tilde{\mathbb{M}}$, the  corresponding row
and column indices and then call the  MATLAB {\tt
  sparse} function, that uses a  compressed sparse column format.
\end{remark} 

%The resulting ratio factor in this case is $\approx 10^5$.

In the last test, we experimentally study the growth order of the computational effort needed to form $\tilde{\mathbb{M}}$, and we highlight which parts of the code mainly contributes to this effort.

\begin{figure}
\begin{center}
 \includegraphics[width=0.6\textwidth] {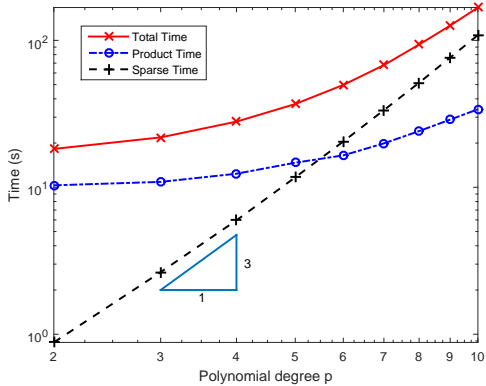} 
\end{center}
\caption{ Time for mass matrix formation with the WQ approach proposed in this paper. In this case $d=3$, $n_{el}=40^3$. Reference slope is $p^3$. 
 Along with the total time, we show the time spent by the product \eqref{eq:sumfactorization} and by the function {\tt sparse}, which represent the single most relevant computational efforts of our code. Other timings, which become negligible for large $p$, are not shown. 
% COMMENTO SUI TEMPI TRASCURATI \B 
} 
\label{figure:3}
\end{figure}

In Figure \ref{figure:3}, we plot in a log-log scale the total
computation time spent by Algorithm \ref{algo_iniz}--\ref{algo_mass} for $40^3$
elements and spline degree up to $10$. We also plot the time spent in
the computations of the matrix-tensor products (i.e., line 5 of Algorithm
\ref{algo_mass}, which is the dominant step  with respect to  the number of FLOPs of
the whole procedure), and the time used by  the MATLAB function {\tt
  sparse}, which is responsible  of allocating the memory for
$\tilde{\mathbb{M}}$ and copying the entries in the sparse matrix data
structure.   These timings were obtained using the profiler of
MATLAB\footnote{ The same timings were also computed using the commands {\tt tic} and {\tt toc}, yielding similar results}.
If we consider the products time, we can see that the its growth
relative to $p$ is significantly milder than what is indicated by the  theoretical FLOP
counting, i.e., $O(N_{\mathsf{DOF}} \, p^4)$. This is probably related to the small dimension of the matrices and tensors involved. 
On the other hand, the times spent by the {\tt sparse} function is
clearly proportional to $p^3$, as highlighted in the plot by a
reference triangle with slope $3$. This is expected, as the number of
nonzero entries of $\tilde{\mathbb{M}}$ is $O(N_{\mathsf{DOF}}\, p^3)$. 
What is surprising is that, for $p> 5 $ the time of the {\tt sparse}
call dominates the total time of the algorithm.  
This indicates that the our approach  is in practice giving the best
possible performance at least for degree high enough, since the  {\tt sparse}
call is unavoidable and well optimised in MATLAB. 

Furthermore, the computing time   depends linearly on
$N_{\mathsf{DOF}}$, as expected, but for brevity we do not show the
results.

\section{Conclusions} 
The proposed algorithm for the formation of isogeometric Galerkin
matrices  is based on three concepts. First, we use a row
loop instead of an element loop. Second, we use WQ that gives
significant savings in quadrature points. Third, we exploit the
tensor-product structure of the B-spline basis functions, 
adopting an optimized  sum-factorization
implementation as in \cite{Antolin:2015}. Our approach
also incorporates an idea  of a previous work: the
numerical computation of univariate quadrature rules as in
\cite{Hughes:2010} and following papers.
In the present work, however, we fix a priori the quadrature points
% in order to have a liner problem to be solved 
so that the weights are given by solving a linear problem,
and we use sum-factorization cycling on rows and not on elements.
The result is a significant gain in performance compared to standard
approaches, for all polynomial degrees but especially for high
degree. For example, in the numerical tests that we present, for $p=6$
the time of formation of a mass matrix is seconds vs hours (comparison
made with GeoPDEs 3.0, which has a well optimised but standard design,
see \cite{de2011geopdes,geopdesv3}), where in our algorithm the
computational time is dominated by the unavoidable MATLAB {\tt sparse}
function call. These results pave the way to the practical use of high-degree
$k$-refinement.  Moreover they relight the interest for a
comparison between Galerkin and  collocation formulation, that is nowadays
preferred for high-degree isogeometric simulations, see
\cite{Schillinger:2013}.  Curiously, a  Galerkin formulation with WQ
is closer to collocation, from the
viewpoint of the computational cost and since both do not preserve symmetry,
i.e. the matrices formed from symmetric differential operators  are
not symmetric. However WQ should preserve the other main properties of
Galerkin formulations.

 Our work will continue in three different
directions. We need to develop a full mathematical analysis of this
approach. We will work on a full implementation within
GeoPDEs. Finally, we will develop the proposed approach in the
direction of non-tensor product spaces (T-splines, hierarchical
splines, etc.), where we expect that some significant advantages of
our approach will be maintained. 

\section*{Acknowledgements}
The authors would like to thank Rafael V\'azquez for fruitful discussions on the topic of the paper.
 Francesco Calabr\`o
was partially supported by INdAM, through GNCS research projects.  Giancarlo Sangalli and Mattia Tani
were partially supported by the European Research Council
through the  FP7 ERC Consolidator Grant n.616563 \emph{HIGEOM}, and by the Italian
MIUR through the PRIN  ``Metodologie innovative nella modellistica
differenziale numerica''.  This support is gratefully acknowledged.

\end{document}